\documentstyle[12pt]{article}
\def\={\!=\!}

\newcommand{\spacenine}{{\ \ \ \ \ \ \ \ \ }}

\newcommand{\ep}{{\epsilon}}

\newcommand{\tilder}{{\tilde{r}}}

\newcommand{\tildew}{{\tilde{w}}}

\newcommand{\E}{{\rm \bf E}}

\newcommand{\argmax}{{\rm argmax}}

\newcommand{\esc}{{\rm esc}}

\def\={\!=\!}

\def\|{{\Vert}}

\def \esc {\mbox {esc}}

\def\ep{\epsilon}

\def\noi{\noindent}

\def\skip{\vspace{2mm}}
\newcount\u
\def\Remark{\skip\noi{{\bf{Remark \number\u.}}} \advance\u by 1}
\date{}
   \title{A Proof of the Vieille result using a kind of discount factor}
\author {Robert Samuel Simon}

\begin{document}
\maketitle

\thispagestyle{empty}

\vfill     
   
\setcounter{page}{-1}

\newpage

{\bf Abstract:}  \vskip.2cm 

We give an alternative  proof that  
  every  two-person 
 non-zero-sum absorbing positive recursive  
 stochastic game with finitely many states 
 has approximate equilibria, a result proven 
 by Nicolas Vieille. 
  Our proof uses    
 a  state specific  discount factor which is similar to the 
 conventional discount factor only when there is only one non-absorbing 
 state.  Additionally  we show 
 that if the players engage in  time homogeneous Markovian 
 behavior relative to some finite state space of size $n$ then 
 for the existence of an $\ep$-equilibrium it 
 suffices that   one-stage deviation brings no more than 
 an $\ep^3/(nM)$ gain to a player,
 where  $M$ is a bound on the maximal difference between 
 any two payoffs.    \vskip1cm

{\bf Key words}: Stochastic Games, Markov Chains

\vfill     
   
     This  paper would not have been possible  without the very 
 generous help  of 
 Eilon Solan; furthermore some of the paper    
  was written by him.    
   This research  was supported by 
     the Center for Rationality and Interactive Decision Theory (Jerusalem),  
    the Department of Mathematics of the Hebrew University (Jerusalem),    
 the Edmund Landau Center for Research in Mathematical
Analysis (Jerusalem)  
       sponsored by the Minerva Foundation (Germany),  the German Science 
 Foundation (Deutsche Forschungsgemeinschaft), and the Institute of 
Mathematical Stochastics (Goettingen).

 \thispagestyle{empty}

    \newpage


\section{Introduction}

A two-player stochastic game is played in stages.
At every stage the game is in some state of the world.
Both players are informed of the whole history,
including the current state, and based on this information 
 they choose simultaneously a pair
of actions.
The current state and the pair of actions chosen determine both
a stage payoff for each of the players
and a probability distribution according to which a new state 
is chosen.

For any $\epsilon \geq 0$, an  
{\em $\epsilon$-equilibrium} in a game 
is a set of strategies, one for each player, 
 such that no player can gain in payoff by more than $\epsilon$ by choosing 
  a different strategy, given 
  that all the other players 
  do not change their strategies.  A game has approximate equilibria if 
 for every positive $\epsilon>0$ it has an $\epsilon$-equilibrium. 
 The {\em value} of a zero-sum game, should one exist, is the 
 unique cluster point of the $\ep$-equilibrium expected payoffs (for 
 the first player) as $\ep$ 
 goes to zero.     
     The un-discounted  payoff of a player in a stochatic game with infinitely 
      many stages, when defined, 
      is a limit  as the number of stages goes to infinity   
    of the average summed over the stages of the player's  expected  payoffs.
   Unless specified, the payoffs of a stochastic game are undiscounted.

Shapley (1953) presented the model of stochastic games,
and proved that a discounted  zero-sum games
always have a value obtainable with 
 stationary optimal strategies.
This result was generalized for equilibria in  $n$-player non-zero-sum 
 discounted  games by Fink (1964).

An absorbing state is such that the play never leaves this state
once it is reached.   
Kohlberg (1974) proved that every two-player zero-sum stochastic game with 
only one non-absorbing state has a value. 
Based on the work of Bewley and Kohlberg (1976),
 Mertens and Neyman (1981) generalized this result, and proved
that every zero-sum stochastic game has a value.

A stochastic game is {\em recursive} if the stage payoff at all non-absorbing 
 states is zero, no matter what the players do.  A recursive stochastic game 
 is {\em positive} recursive if there is a player who receives at all 
 absorbing states only positive 
payoffs.  A positive recursive 
 stochastic game is {\em absorbing} 
  if the player who receives these positive  payoffs 
 can force the play toward absorption.

Existence of  approximate equilibria  in two-player non-zero-sum
stochastic games with only one 
non-absorbing state was proven by   Thuijsman and Vrieze (1989).
In their proof  Thuijsman and Vrieze considered a sequence of stationary 
 equilibria of the discounted game as the discount factor tends to 1,
and they constructed different types of $\ep$-equilibrium strategies
according to various properties of the sequence.

Vieille (2000a) showed that for approximate equilibria to exist in  every
two-player non-zero-sum stochastic game with finitely many 
 states 
it is sufficient to prove this  for the sub-class of
absorbing  positive recursive games.
Furthermore Vieille (2000b, 2000c) proved that indeed all games in 
 this subclass have approximate equilibria.

In the present paper, we provide an alternative proof of the Vieille 
 result for  absorbing  positive recursive games.
 The primarily difference between our proof and Vieille's lies 
 in the use of a kind of discount factor 
 rather than Vieille's undiscounted evaluation. 
This discount factor is state specific and is similar to the conventional
 discount factor only when there is only  one non-absorbing state.  
We were inspired by the  Thuijsman and Vrieze article  
 and their confidence  that  their 
 ideas could deliver the same result for 
 finitely many states.   Our goal was to confirm their 
 optimism by demonstrating  
 the great versatility of the discounting concept.

In positive recursive games, discount factors for the player receiving positive
 absorbing payoffs  
 persuade him to make  moves that push the game 
toward absorption.   Let us call this player the {\em second} player.
The  serious problem  with generalizing  
 the Thuijsman and Vrieze approach directly   
  is that 
      the usual discounted evaluation does not discriminate 
       between the time spent at the state at which a 
       decision is made and the other 
        states that might follow this decision. 
 As long as the second player  at a given state 
 chooses between two moves that do not involve returning to that state, 
 his  evaluation of those moves in an appropriate discounted game  
 \underline {should} be based upon his 
 undiscounted evaluation. 
  Play that never returns to  
      this  state  before  
       absorption but   visits 
       other states arbitrarily 
       many times receives no discount whereas 
     play that re-visits the initial state $n$ times receives  
      a $(1-\delta)^n$ discount, 
regardless of its visits to other states. 
  
 We see no way  to generalize our  proof   to 
 three player games   (and  
 it appears highly unlikely). 
    On the other hand, we can not dismiss the possibility; (see also 
 Solan, 1999, where discounted evaluations were used to understand 
 some three player undiscounted stochastic games).      
 If the compactification of a   strategy 
 space  creates discontinuities in the undiscounted 
  payoffs   a discounted 
 evaluation   may  handle the points of discontinuity  
 successfully.   A false impression that discounting is 
useless to understanding the undiscounted game 
 may result from a lack of knowledge of 
 how to {\em turn off} the discount where one is sufficiently far from 
 the points of discontinuity.     As we will see below, knowing when to 
 turn off the discount is central to our approach.

The secondary difference between our proof and Vieille's is that  the 
 mathematics we use  is entirely   elementary.  No deep theorems 
 of mathematics are required; for example, there is no use of the theory of  
  semi-algebraic functions. What we need from the theory of Markov 
 chains is very elementary and proved entirely 
 in this paper. Due to our discounting approach we work  
  with 
 taboo probabilities rather than the directed graphs  perspective 
 of Freidlin and Wentzell, (1984).

The only theorem we quote instead of proving is Doob's  submartingale 
 inequality, a generalization of Kolmogorov's inequality and also 
 an easy theorem to prove. Applying the inequality, we show 
 that if the players engage in  time homogeneous Markovian 
 behavior relative to some finite state space of size $n$ then 
 for the existence of an $\ep$-equilibrium it 
 suffices that   one-stage deviation brings no more than 
 an $\ep^3/nM$ gain to a player,
 where  $M$ is a bound on the maximal difference between 
 any two payoffs.   
   \vskip.2cm 

{\bf Countably many states}
\vskip.2cm 
 
We developed our  unorthodox approach to stochastic games   
   with the hope that it would deliver   
 approximate equilibrium existence for
 all two-person non-zero-sum stochastic games 
 with countably many states. 
  We have failed  in this attempt.

The main problem  is that our approach (and that of Vieille) rests 
 ultimately on the pideon-hole principle.
   If the expected number of 
 visits to every non-absorbing 
state is finite  then 
 with probability one  an absorbing state is reached.   This does not hold 
 if there are infinitely many non-absorbing states.

  In general, what is the difficulty 
 in proving approximate equilibrium existence 
 for non-zero-sum two-person  stochastic games with countably many 
 states?   Several important positive results need to be mentioned.
Maitra and Sudderth (1991) proved that all zero-sum stochastic 
 games with countably many states have values.
 In a game of perfect information, the players take turns making their 
 moves and each player knows the previous moves of the other players; 
 the classic example is that of chess.
 A {\em Blackwell game} is identical in transition structure 
 to a stochastic game, but the payoffs are determined by 
   a  function Borel measurable with 
 respect to the histories of play.
Martin (1975) proved  that all zero-sum  Blackwell games of perfect 
 information have values, and 
  Mertens and Neyman (in Mertens 1987)
 extended Martin's result to  non-zero-sum games with finitely many players. 
  Using his  result for games of perfect information, Martin (1998) 
 proved that all 
 zero-sum Blackwell games have values.

The differences between 
  non-zero-sum stochastic games (with simultaneous moves) and either 
 non-zero-sum Blackwell games of perfect information 
 or zero-sum Blackwell games with simultaneous moves 
  are formidable.     The probability of  
  absorption at a stage  in a stochastic game 
  can be also a  minimal bound on that  stage's 
  deviation from pure equilibria; (for example see the ``Big Match" in 
Blackwell and Ferguson, 1968).  With the $\epsilon$-equilibria of many 
 games, including the absorbing positive recursive variety, 
 while absorption must become a near certainty 
 the  culmulative opportunity to exploit  
 deviations must  not exceed  $\epsilon$.   Therefore
 one needs  that stage for stage approximate equilibria can 
 translate  to cumulative approximate equilibria. 
 In zero-sum games this is not so problematic 
  because  the gains to  one player from deviation  
equal the losses to the other player.
  But with two-person  non-zero-sum games, one must consider functions 
 with values in ${\bf R}^2$; the  
   potential independence of 
 the two values and  need for a cooperative solution frustrate 
 attempts to generalize the  approachs that were 
 successful with zero-sum games.   On the other hand if the moves are 
 made simultaneously how does one know the other player is adhering 
 to a cooperative agreement? 
   So far the  main answer   
    has been to 
 request from  each player Markovian behavior, 
  accompanied by statistical 
 testing and punishment by the other player 
 in the event of significant  statistical deviation. 
  With this approach, it is necessary that the probability that   
   an honest player will be punished unjustly can be made arbitrarily 
 small. 
 As we will demonstrate with the following proposition and counter-example 
 to a variation on this proposition, such a control process is  unlikely   
 in general for Markovian behavior that is carried out essentially on 
 a  countable  state space.

If $S$ is a finite or countable set let $\Delta (S)$ stand for the space 
 of probability 
 distributions on $S$. A {\em Markov chain} is defined by
 a finite or countable state space $S$ and 
 for every $s\in S$ and stage $i\geq 0$ a probability distribution 
 $p^s_i\in \Delta (S)$ governing the distribution on the states at the 
 $i+1$st stage, given that $s$ is the state on the $i$th stage. It is 
 {\em time homogeneous} if $p^s_i$ is indendendent of the $i$. 
\vskip.2cm

{\bf Proposition 4.2:} 
 Let $X$ be a finite  space. For 
 every $x\in X$ let  $Y_x$ be a finite space, with $Y:= \cup_{x\in X}Y_x$. 
  (In the context of stochastic games, $X$ will be the state space and 
$Y_x$ will  be the set of 
 moves that a player has at the state $x\in X$.)  There are
  probability transitions
 $(p^x\in \Delta (Y_x)\ | \ x\in X)$ from $X$ to $Y$ 
 and  there are probability transitions 
 $(p^y\in \Delta (X) \ | \ y\in Y)$
 from $Y$ to $X$, so that    
  for every starting point $x_0\in X$  a time homogeneous 
 Markov chain on $X\cup Y$ is defined. On the even  stages $i=0,2,4,\dots$ 
the process is in 
 $X$ and on the odd stages the process is in $Y$. 
 Let  there be an  evaluation function $v:X\cup Y \rightarrow 
 {\bf R}$       that is harmonic  with respect to   
  the transitions (meaning 
 that  a  martingale is formed). Let 
 $M>0$ be a  uniform 
 bound for the  maximal difference between all values of  $v$.
For every pair 
  $x\in X$ and  $y\in Y_x$ such that $y$ is 
 reached from $x$ with positive probability (according to $p^x$)
 the 
 difference between
 $v(y)$ and $v (x)$ is no more than $\delta > 0$.\newline  
{\bf  Conclusion:}  If $|X|=n$, $\epsilon < 1/2$, and       
  $ \delta  \leq  \epsilon^3 /Mn$ 
  then  
 the probability that there exists an $l$ with 
 $\sum_{i=0,2,\dots}^l (v(y_{i+1})-v (x_i))\geq \epsilon$ 
  does not exceed $\epsilon$. 
\vskip.2cm 
The complexity   of the $Y_x$ play no role in the proof of 
 Proposition 4.2, and therefore it could have many 
  generalizations corresponding to variations in the structure of the $Y_x$.

  To emphasize the importance of the finite number $|X|$, 
 the following is  a counter-example to Proposition 4.2 if   
 we assume that the bound for $\delta$ is   
independent of the cardinality of $X$. Furthermore, 
 if we consider processes that are 
 not time-homogeneous, it does not 
 help  if for every stage  the sum over the states of the maximal  
 differences   add up 
 to no more than  $\delta$.    

Consider a random walk 
 on $n+1$ positions such that at the left end (at position $0$)
 the player receives an absorbing payoff of 
 $0$ and on the right end (at position $n$) 
an absorbing payoff of $1$.  The space $X$ is the $n+1$ positions and 
 for every $x\in X$ the two-set $Y_x$ consists of the two directions 
 ``left" and ``right". Given any small $\delta>0$, 
 one can make $n$ large enough so that at every stage the change in expected 
 payoff does not exceed $\delta$.  Now reformulate the randon walk 
 so that at the  $k$th  stage of play 
 there is no motion at any $i$ position with  $i\not=k$ (mod $n-1$),
 but at the $k'=k$ (mod $n-1$) 
position there is an equal $1/2$ probability of moving 
  either to the position $k'-1$ or to $k'+1$.  At each stage 
the sum over the states of the differences in expected payoffs 
  remains no more than $\delta$, and yet we are no closer to satisfying 
 the conclusion of 
 the proposition.  (With $n$ even and  
starting in the middle position 
 with an expected payoff of $1/2$, for every small positive $\epsilon$
 with   probability close to $1/2$ there will be motion to a position 
 with an expected payoff of at least $1/2 + 2\epsilon$.)

  We expect no proof of approximate equilibrium existence for all 
non-zero-sum  stochastic games with countable state spaces 
without  a radically different approach. 
 If a proof for countably many   states can be  
 found, its application to finite state truncations of the countable 
 state  game 
   would  provide approximate  
 equilibria of the finite state games such that the average number of 
 stages before absorption would not  
 explode with the increase in 
 the finite number of these states.  In the 
  proof below 
 for a fixed $\ep$  there is   no lower bound determined by  the 
 number of states on the rate   for which an absorbing 
 state is reached.  
  Indeed, because such a proof would imply the existence 
 of yet another  alternative proof for finitely many states with   dramatic 
 absorption rate properties, we  
 suspect  that there is a counter-example. 
 Furthermore, it is
 possible that the complexities from countably many states 
 involved in a two-player counter-example  
  could be mimicked  by 
 the introduction of more players in a stochastic  game with 
finitely many states, yielding a counter-example to approximate 
 equilibria in this context as well.

 We suspect that approximate 
equilibrium existence for a broad class of two-person 
stochastic games played on countable state spaces must rest on a 
fundamental  assumption: 
 that  there is a uniform bound on 
 the number of states possible on any 
 given stage of play.           
  With a finite number of such positions, it is still not clear
 how appropriate Markovian   
  should be found.  
 Even with only one non-absorbing position, the possible infinite variations, 
 including  the number of moves for each player and the order in which similar 
 ``types"  may appear, 
  make the problem formidable.    At least the generalization of 
 Lemma  4.1 to Markov chains that are not time homogeneous will be 
 necessary.  
Another reason to  present our alternative proof of the Vieille result 
 is the hope that it will be relevant 
 to this case, which we call the case of {\em finitely many positions}.
   If for each non-absorbing  position 
one could  find an appropriate  common  identity to an infinite 
 sub-sequence of states occuring in that position, then the pideon hole 
 principle could be applied successfully.
  Throughout this paper, we 
 comment on the case of finitely many positions.

 \vskip.2cm 

{\bf Organization}
\vskip.2cm

To execute our  proof  efficiently,
 we will assume that Player One 
 has the ability to send signals to Player Two that are independent of  
 the transitions in the games.   The easiest way to formalize this 
 property is to assume that every move of Player One at a non-absorbing state 
 is paired with another move at the same state that is its identical 
 copy with respect to the transitions.  Without this assumption, 
 the proof is  
  formality  more involved,  less elegant, however 
 essentially equivalent.  In the  section following 
 the conclusion of the main proof, we 
 prove the result    
 without this signaling assumption.

The argument and the paper  are organized as follows. 

Section 2 introduces the model of absorbing positive 
 recursive stochastic games and the basic concepts of Markov 
 chains. Additionally we 
 introduce an important concept with regard to the movement between 
 states, called {\em taboo} probabilities. 
 A taboo probability is the probability
 that one moves from an initial state to some  set of target states 
without travelling through 
 some second set of ``forbidden" states.

Section 3 gives proofs of all 
 the needed lemmas on Markov chains.   
   The  most central lemma is Lemma 3.2; it 
  states that when motions at a multitude 
 of states  are removed whose frequencies are only a small fraction 
 of  the total  motion
  toward a fixed state     then 
  the flow continues  toward this fixed state 
with about  the same or greater tendency.

Section 4  contains a proof of Proposition 4.2, which 
 also establishes  general  sufficient 
 conditions for the existence of approximate equilibria.   
We create new states from our old states,
 which we call {\em situations}; at most three situations are created 
 from each original state.  The method of creating 
 the situations we call {\em polarization,} introduced  in 
 Section 3.  Except
 for the rare possibility of punishment, our 
 behavior strategies will be stationary on the situations. 
 Section 4 concludes with Theorem 1, a demonstration of sufficient 
 conditions for approximate equilibrium existence in our games.

 In Section 5 we introduce the state specific discounted evaluation 
 for the second player. 
  We define the discounted evaluation 
  such that the discounting rates are adjusted for   
states sufficiently close together, according to a  metric determined 
 by the strategies. 
  We select a quantity $\overline\ep$ much smaller than $\ep$, and 
 define the discounted evaluation so
  that moves with more than an $\overline\ep$ probability
 of non-return to the state are evaluated in an undiscounted way 
 and moves with a $\gamma$ probability of no return with 
 $\gamma< \overline \ep$ are evaluated as if 
 their probability of no return was $\gamma /\overline \ep$. Our 
 choice for $\overline\ep$ is guided by  Proposition 4.2.

 A serious problem with the state specific discounted evaluation 
 is that  the motivations of the 
 second player at one state can be very different from that at another state. 
 Essentially the second player becomes a multitude of players, one 
 for each state.   This allows for the second player  at some
 states  to prefer moves that result in too slow a motion toward 
 absorption and therefore also discounted evaluations below 
 the zero-sum value.  To avoid this problem, in Section 2 we  
 define a new correspondence, called the ``jump" correspondence,  based upon  
 stationary strategies  optimal in the conventionally discounted game.  The 
 use of the jump  
correspondence  by the second player 
 results in fast absorption.  The   ``best-reply" correspondence 
 of the second player is a combination of the jump correspondence with 
 a maximization of the state specific 
discounted evaluation -- when the discounted evaluation 
 is too low, the jump correspondence is activated. 
 For the first player, the undiscounted evaluation is used to define her 
 ``best-reply" correspondence.  With the 
  ``best-reply" correspondences for both player defined,
 we demonstrate two important properties.   
Lemma 5.4 shows that at a fixed point
 the jump correspondence of the second player has 
 only  very limited influence on the play. Lemma 5.5 contains 
 the key argument to our entire approach;
 it is used repeatedly to solve the most difficult problems. 
  It shows that if there is a  meaningful discrepancy 
 between   the discounted and undiscounted evaluations for the second 
 player  then the second player  seeks primarily 
 motion with the fastest  absorption rate.

The synthesis of the previous sections lies in Section 6. Theorem 
 2 proves that the conditions of Theorem 1 are always satisfied -- implying 
  the existence of approximate equilibria. 
  Here we consider 
 sets such that   a significant proportion 
 of all the motion leaving these  sets are from   
 Player Two  moves with payoffs for Player Two 
  significantly below the set-average payoff. 
 Fixing any such state in a  set where such moves take place, 
   we look at what happens 
 when Player One stops playing all moves performed with   
  frequencies small compared to the motion toward this special state.
 The result, for which Player One is indifferent,  involves 
almost exclusively the use of similar such moves by Player Two such that 
 the players can travel between these moves without the danger that 
 along the way  Player Two prefers to provoke punishment over performing 
 one of these moves.  Ultimately we 
 show that there is a convex combination of such moves that all yield the same 
 payoff for Player Two and for which Player One is approximately indifferent.

In Section 7 we consider the problem of signaling, as described above; and 
 in Section 8 we conclude in more detail 
 with the problem of countably many states.

\section{Preliminaries}
    \subsection{ The Model}

Let $\cal S$ be  the set of   states;
 $\cal A $ is the subset of absorbing states and 
 ${\cal N}= {\cal S}\backslash {\cal A}$ 
is the subset of non-absorbing states.

       For every $s\in {\cal S}$,   
            $A_1^s$ and 
      $A_2^s$ are the moves (pure actions) 
 of the first and second players,  
              respectively, at the state $s$. Without loss 
 of generality, we assume that $|A^s_i|=1$ for every $s\in {\cal A}$ and 
 $i=1,2$. 
           Let                  $r^1:{\cal A}\rightarrow [-1/2,1/2]$ 
                 and $r^2:{\cal A}\rightarrow [\omega,1]$ 
                 be the first and second players' 
                  evaluations on absorbing states, respectively, with   
  $0<\omega < 1$.
       Let $m$ be the maximal number of moves for either player at 
       any non-absorbing state,
  meaning $m=\max_{s\in {\cal N}} 
                (|A_1^s|, |A_2^s|)$.
                
Let $p(t|s; a,b)$ be the probability 
               of moving from $s$ to $t$ when $a\in A^s_1$ and $b\in 
                A^s_2$ are played.  Let 
 $\rho$ be  
defined by $\rho := 
                 \min ( 
  p(t|s;a,b) \ |\ s,t\in {\cal S}
 \ p(t|s;a,b)>0)$, the minimal 
 non-zero transition probability. Notice that in the case of finitely 
many positions  
  one has such a positive quantity for each stage.  More relevant, however, 
 would be a sequence $\rho_i$ of positive quantities such that the 
 series $\rho_i$ is divergent but sums  toward infinity 
 much slower  than any 
 divergent series of positive transition probabilities.  Such a series 
 is possible if there is  a  uniform bound on the number of moves.     
 Additionally the discount factor must be adjusted to this series, 
 (possibly with the discount factor equaling $1-\delta \rho_i$ if 
 there is only one non-absorbing state).

Let $X:= \prod_{s\in {\cal N}} \Delta (A^s_1)$ and 
$Y:= \prod_{s\in {\cal N}} \Delta (A^s_2)$
 be the spaces of stationary strategies of the  
  players, with $X^s:=\Delta (A^s_1)$ and 
  $Y^s:=\Delta (A^s_2)$. 
            For $a\in A^s_1$, $b\in A^s_2$, 
            $x^s\in X^s$ and $y^s\in Y^s$ we define $p(t|s;a,y^s)$,
              $p(t|s;x^s,b)$ and $p(t|s;x^s,y^s)$ in 
               the appropriate linear or bi-linear way.
For any $s\in {\cal N}$, $x^s\in X^s$ and $a\in A_1^s$, 
 the quantity $x^s_a$ will 
 stand for the probability, as determined by $x^s$, 
 that the move $a$ is used. The same applies for 
 $b\in A_2^s$, $y^s\in Y^s$ and $y^s_b$. 
Define a pair $(x,y)\in X\times Y$ to be {\em absorbing} if from every start
  with probability one an absorbing state is reached.

We will say that two positive quantities $a$ and $b$ are 
 different by no more than a   
  {\em factor} of positive 
            $\gamma<1$ if  $a\geq b (1-\gamma)$ and 
 $b\geq a(1-\gamma)$.

\subsection {Histories, Strategies, Equilibria}

For every stage $i\geq 0$ and $s\in {\cal S}$  define ${\cal H}^s_i: = 
\{ (s_0, a_0, b_0), (s_1, a_1, b_1), \dots ,\newline  
(s_{i-1}, a_{i-1},b_{i-1}), s_i=s\
|\ \forall\  0\leq k <i\ 
\ a_k\in A^{s_k}_1, b_k\in A^{s_k}_2,\  p(s_{k+1}| s_k; a_k,b_k) 
 >0\}$, 
 with ${\cal H}^s_0= \{ s\}$ for all $s\in {\cal S}$. 
Define ${\cal H}^s := \cup_{i=1}^{\infty} {\cal H}^s_i$,  
 ${\cal H}_i:= \cup _{s\in S} {\cal H}_i^s$, 
${\cal H}:=   \cup_{i=0}{\cal H}_i$, and 
 $\tilde   {\cal H}:= \{ (s_0, a_0, b_0), (s_1, a_1, b_1),
 \dots \
|\ 
 \forall i\geq 0 $ the truncation up to $s_i$ belongs 
 to ${\cal H}_i^{s_i} \}$, the 
 set of infinite sequences.  

A strategy of  Player $j=1,2$ is a set of maps $\sigma_j=(\sigma^s_j\ | \ 
 s\in {\cal N})$ with $\sigma^s_j$ a map from 
 ${\cal H}^s$ to $ \Delta (A_j^s)$  for all $s\in {\cal N}$.

With Blackwell games, a more general class than stochastic games,  
we assume that a player's evaluation on $\tilde {\cal H}$ 
  is a function 
 that is measurable with respect to the Borel subsets of $\tilde {\cal H}$, 
 the sigma algebra induced by the subsets of ${\cal H}_i$ for all $i\geq 0$.  
In case that a stochastic  
game is recursive, for every member of $\tilde   {\cal H}$ it easy 
 to define an evaluation for both players.  Either the infinite 
 sequence reaches an absorbing state and
 the players receive the corresponding
 absorbing 
 payoffs,
 or it never reaches an absorbing state and both players receive 
 a payoff of zero.    

For every initial  state $s$  and 
every pair of strategies $\sigma_1,\sigma_2$ for both players 
 a distribution is induced on $\tilde 
{\cal H}$ in a natural way, resulting in   
  two evaluations ${\cal V}^s_j(\sigma_1, \sigma_2)$ for Player $j=1,2$ of 
 the expected values of the $r^j$ on $\tilde {\cal  H}$. 
  An $\epsilon$-equilibrium 
 is a pair $\sigma_1,\sigma_2$ such that for all $s\in S$ and 
 alternative strategies 
 $\tilde   \sigma_1$ and $\tilde   \sigma_2$ it holds that 
 ${\cal V}^s_1(\tilde   \sigma_1, \sigma_2)
\leq {\cal V}^s_1(\sigma_1, \sigma_2)+\epsilon$ and 
  ${\cal V}^s_2( \sigma_1, \tilde   
\sigma_2)\leq {\cal  V}^s_2(\sigma_1, \sigma_2)+\epsilon$. 
With absorbing positive recursive games and positive $\omega$  the 
 lowest Player Two absorbing payoff  we get the additional 
 property that there exists an $N>0$ such that with probability 
 at least $1-{2\epsilon \over \omega}$ the game has reached 
 an absorbing state before the stage $N$.

\subsection{Jump Function}
 For any  positive 
 real number $0<\alpha<1$   
let ${\cal G}^{\alpha }$ be the conventionally defined 
 discounted  zero-sum game played 
 against Player Two   
  such that a visit to any state is discounted 
 according to $1-\alpha$, and let ${\cal G}^0$ be 
 the corresponding undiscounted zero-sum game.
 For all  positive  $\alpha$  we define 
 $c ^{\alpha}:{\cal S}\rightarrow {\bf R}$ to be the min-max value for 
 Player Two in 
 the zero-sum game ${\cal G}^{\alpha}$,   
 with $c^{\alpha}(s)=r^2(s)$ for all $s\in {\cal A}$.   
 Because the game is positive recursive 
 the $c^{\alpha}$ are  monotonically non-decreasing 
and due to Mertens and Neyman (1981) the point-wise limit is  
 the undiscounted value of the game ${\cal G}^0$, though for 
this class of games there is an elementary proof.  Player Two chooses 
 a stationary optimal strategy of ${\cal G}^{\alpha}$ for an $\alpha>0$ 
 sufficiently small so that $c^{\alpha}$ is within $\ep$ of 
 its point-wise limit and at stage $i$ Player One chooses one of her optimal 
 strategies in the game ${\cal G}^{\alpha_i}$ where for every 
 $i\geq 0$
 $c^{\alpha_i}$ is within $\ep/2^{i+2}$ of the point-wise limit and 
  $\alpha_i <\ep/ 2^{i+2}$.

 
For every $x \in X$ and positive $0<\alpha<1$ 
define the {\em jump} function $j^{\alpha}_{x} : {\cal N} \rightarrow R$ by
\[ j^{\alpha} _{x}(s) = (1-\alpha)
 \max_{b \in A^s_2} \sum _{t\in {\cal S}}p^{}(t|s,x,b) \ c^{\alpha}(t) \]
--   the maximal payoff that Player Two can guarantee himself in 
the $1-\alpha$ discounted game   
 by being punished 
 after the next stage if Player One uses $x$ at the present stage.
     If $s$ is an 
 absorbing state, define $j^{\alpha}(s)$ to be $r^2 (s)$ for 
 all $\alpha$. 
For all states  it is clear that 
$j^{\alpha} _{x} \geq c^{\alpha}$,
   with equality when $x$ is 
    an optimal strategy for Player One in the 
 zero-sum game ${\cal G}^{\alpha}$ played against Player Two.  
For every state $s\in {\cal N}$ and $x \in X$ define
\[ J_x^{ \alpha}(s) =
 \argmax_{b \in A^s_2} \sum _{t\in {\cal S}}\ p(t|s,x,b)\  c^{\alpha}(t). \]
Let $n(s)$ denote the state following $s$, in our 
 context a random variable.
If $s$ is not an absorbing state and $b\in J^{\alpha}_x(s)$
 then 
$j^{\alpha}_{x}(s) \leq (1-\alpha) \E^{x}_b j^{\alpha}_{x}(n(s)), $
where $ \E^{x}_b$ is the expectation determined by the move 
 $b$ and the strategy $x^s$.
This makes $j^{\alpha}_x$ a sub-martingale.

For  $i=1,2$ and a state $s\in {\cal S}$
 define $c_i(s)$ to be the  
 value for Player $i$ of the 
zero-sum undiscounted game played against Player $i$ starting at the state 
 $s$.   
For every Player $i$ and every stationary strategy 
 $z$ of Player $k\not=i$  
define the {\em jump} function $j^i_{z} : {\cal S} \rightarrow {\bf R}$ by
$$ j^2 _{z}(s) = 
 \max_{b \in A^s_2} \sum _{t\in {\cal S}} p(t|s,z,b) c_2(t)
 \quad \mbox { or }\quad 
j^1 _{z}(s) = 
 \max_{a \in A^s_1} \sum _{t\in {\cal S}} p(t|s,a,z) c_1(t)$$
--   the maximal payoff that Player $i$ can guarantee himself against 
 $z$ if he is punished on the next stage.   

\subsection {Taboo probabilities}  

  For any time homogeneous 
 Markov chain, a state $s$, and two disjoint sets $A$ and $B$  
 of states  we introduce 
 the ``taboo" probability $P^A(s,B)$ to be the probability, 
 with a start at the state $s$, of reaching 
  the set $B$ before the set $A$ at any stage following the initial stage 
  at $s$.  With $t_C:= \inf \{ n\geq 1\ | \ s_n\in C\}$ 
 $P^A(s,B)$ measures the event that $t_B<\infty$ {\em and} $t_B<t_A$  
  conditioned on   $s_0=s$. 
  If either set is a singleton, we can write its single 
  member   instead  
   of the set. If there is ambiguity concerning which state space or which 
 transitions, we identify them with a subscript. In our context of
  stochastic games and stationary strategies, 
       $P^A_{ x,y}(s,B)$ will be the taboo probability corresponding 
       to the time homogeneous Markov 
       chain generated by $( x,y)\in X \times Y$.

Define a state of a time homogeneous Markov chain to be {\em absorbing} 
 if once this state is reached then the motion remains in this state 
 forever.  The Markov chain is absorbing if for any start  
 with probability one  an absorbing state is reached. 

Before moving toward the proof, we must present some basic notions using 
 the taboo probabilities. These quantities will be defined first 
 for 
time homogeneous Markov chains and then applied 
 to the games.   

For any part $p$ of a transition at a state $s$ or an alternative 
 transition $p$ for that state  
 define $g(p)$ to be the probability that there is no return to 
 $s$ if  $p$ is used at $s$ and
  the transitions remain constant 
 at all other states.  If $p$ was a part of the 
 transition at $s$ then define  $f_p$ to be the  
 frequency with which $p$ is used at the state $s$.
For every choice $(x,y)\in X \times Y$ and pair  $a\in A^s_1$ and  
 $b\in A^s_2$ of moves at the state $s\in {\cal N}$ 
 $g_{x,y}(a,b)$ 
  is the probability that there is no 
  return to $s$
 given                             that  
Player One and Player Two at $s$ play 
the actions $a$ and $b$,  and elsewhere in the future  
 the stationary strategies $(x,y)$.
For a move $b\in A^s_2$ of the second player, define $g^b_{x,y}$ 
 to be $\sum_{a\in A^s_1} x_a^s g(a,b)$, and define 
 $g^a_{x,y}$ for all $a\in A^s_1$ correspondingly.

Define the absorption rate $a(s)$ of a state $s$
  to be the probability that after any visit to this state 
 there is no return to this 
 state, meaning  that the absorption rate is the expected value 
 of the function $g$. 
For the game  
  the {\em absorption rate} $ a_{x,y}(s)$ of a state $s$ is 
  $\sum_{a\in A^s_1, \ b \in A^s_2} x^s_a y^s_b g^{s}_{x,y}(a,b). $
Given that $(x,y)$ is absorbing 
   $a_{x,y} (s)$ would be the  taboo probability 
     $P_{x,y} ^s(s,{\cal A})$. 

For any  part $p$ of the transition at a state $s$  
 define $\nu (p)$ to be the 
 probability that at the last visit to $s$ the part 
$p$ was used, or equivalently $\nu(p)=f_pg(p)/a(s)$.
We call this the {\em importance} of $p$. 
For a pair of moves $a\in A_1^s$ and 
 $b\in A_2^s$ at $s\in {\cal N}$ 
and  stationary strategies $(x,y)$  the importance  
  $\nu^s_{x,y} (a,b)$ is 
  $x^s_a y^s_b g_{x,y}(a,b)/ a_{x,y}(s)$. 
  For any move $a\in A_1^s$ define 
 $\nu^a_{x,y}$ to be $\sum _{b\in A_2}  \nu_{x,y}(a,b)=
 x_a g_{x,y}^a/a_{x,y}(s)$ 
 and 
  for any move $b\in A_2^s$ define 
 $\nu^b_{x,y}$ in the same way.

For any distinct pair $s,t$ of states define 
 $\esc (t,s)$ to be the probability of never reaching $s$ with a start 
 at $t$. ($\esc$ stands for ``escape".)  
For the game we have  
   $g^{b}_{ x,y}= \sum_{t\in {\cal S} } 
     p (t|s; x,b) esc_{ x,y} (t,s)$. 
          (If $(x,y)$ is absorbing, ${esc}_{ x,y}(t,s)$ is
 $P^s_{x,y}(t,{\cal A})$ and is  
           different from $P^{s,t} _{ x,y}
 (t, {\cal A})$, the probability 
            of absorbing before returning to either $s$ or $t$).

          For distinct states $s$ and $t$  let  $\mu(s,t)$ be 
               ${esc}(s,t)+{esc}(t,s)$, 
                and otherwise let $\mu(s,s)=0$.       
$\mu$ is a metric on the state space.  
       Recognize $1-esc (t,s)$ as the probability of moving from $t$ to 
        $s$, and for mutually distinct $u,v,w$ we have 
          $1-esc(u,w)\geq (1-esc (u,v)) (1-esc (v,w))\geq 1-esc (u,v)-
          esc (v,w)$.

Given that the Markov chain is absorbing with $ A$ the set 
 of absorbing states,  
the following relations for states $s\not= t$ are easy to verify:
$$ \esc(s,t) = 
{P^{\{ s,t\}} (s,A) \over P^s(s,t) + P^{\{ s,t\} } (s,A)} 
={P^{\{ s,t\}} (s,A) \over 1-P^{A\cup \{ t\}}(s,s)} 
\spacenine \spacenine \hfill (1) $$
$$ a(s) = P^s(s,t) esc(t,s) + P^{\{ s,t\} } (s,A) \spacenine \spacenine 
\spacenine \spacenine \spacenine \hfill (2)$$
which imply $P^s(s,t) \mu (s,t) \leq a(s)\leq  
\mu (s,t)$ and $a(t) P^s(s,t) \leq a(s)\quad \quad \quad  (3)$.

  For all these quantities and following ones,   
     we can drop the subscripts and superscripts if there is no ambiguity.

\subsection {Evaluations}

 We had extended 
 the values $r^i:{\cal A}\rightarrow {\bf R}$ on the absorbing 
 states to  functions $r^i$  on all paths in $\tilde {\cal H}$.
For any stationary strategies $(x,y)$ and players $i=1,2$ extend the 
 definition of $r^i$ again to a harmonic function 
 $r^i_{x,y}:{\cal S} \rightarrow {\cal R}$ with  
$r^i_{x,y}(s)$ equal to the expected value of $r^i$ 
 on $\tilde {\cal H}$ 
 as determined by $(x,y)$. 

For any  harmonic function $r$ on $S$, and $p$ a part  of or 
 an alternative to  the
 transition from a state $s$, define $v^r(p)$ to be the expected value 
 of $r$ conditioned on the use of $p$ and no return to the state $s$, 
 with $v^r(p)$ defined to be $r(s)$ if there is return to $s$ with certainty. 
  If the Markov chain is absorbing and 
 $g(p)>0$ then  
 $v^r(p)$ would be the new harmonic function value for $s$ if 
 the transition 
 from $s$ were replaced by $p$.    
 For every pair of moves $a\in A_1^s$ and $b\in A_2^s$  
   $v^{i}_{x,y}(a,b)$ is defined to be $v^{r^i_{x,y}}$ of 
 the part of the transition defined by the pair $(a,b)$ of moves. 
 Likewise define $v^i_{x,y}(a)$ and $v^i_{x,y}(b)$ with respect 
 to the pairs $(a,y\in Y^s)$ and $(x\in X^s ,b)$, respectively.  
  If $(x,y)$ is absorbing we have the relation 
 $$r^i_{x,y}(s)= {\sum _{a,b} x^s_a y^s_b v^i_{x,y}(a,b) g^s_{x,y}(a,b) \over 
   a_{x,y}(s)}= 
 \sum_{a,b}\nu_{x,y}(a,b) v^i_{x,y}(a,b).$$
 For $a\in A^s_1$ we have  \newline $v^i_{x,y}(a):= 
\sum_{b\in A_2^s}  y^s_b v^i_{x,y}(a,b) g_{x,y}(a,b)/ 
g^a_{x,y}=\sum_b \nu _{x,y} (a,b) v_{x,y}^i (a,b) / \nu _{x,y}^a$\newline 
  and 
for $b\in A^s_2$   we have  \newline $v^i_{x,y}(b):= 
\sum_{a\in A_1^s}  x^s_a v^i_{x,y}(a,b) g_{x,y}(a,b)/ 
g^b_{x,y}=  \sum_a \nu _{x,y} (a,b) v_{x,y}^i (a,b) / \nu _{x,y}^b$,\newline 
 with both quantities 
 $r^i(s)$ when the quotient is not well defined.  

 For any  harmonic function $r$ on $S$, and $p$, a part of or an  
 alternative to  the 
 transition from a state $s$, define  $w^r(p)$ to be  the expected value 
 of $r$ on the following stage according to 
 the   one-time 
 use of $p$ on that stage.
  We have $w^r(p) = g(p) v^r(p)+ (1-g(p)) r(s)$. 
For any pair of moves $a\in A_1^s$ and $b\in A_2^s$ at $s\in {\cal N}$    
  and $i=1,2$  
 $w^{i}_{x,y}(a,b)$ is the expected  value
 of $r_{x,y}^i$  
 on the next stage  if the players use 
 the pair $a$ and $b$ on the 
 present stage at $s$.
  For all $b\in A_2^s$ 
 define $w^i_{x,y}(b):= 
 \sum_{a\in A^s_1} x^s_a w^i_{x,y}(a,b)$ and 
 for all $a\in A_1^s$ define 
  $w^i_{x,y}(a):= 
 \sum_{b\in A^s_2} y^s_a w^i_{x,y}(a,b)$. 
\vskip.2cm 

The following is a  central lemma concerning the changes in a harmonic 
 function. 
\vskip.2cm

{\bf Lemma 2.1:}  Let $S$ be the finite state space of an absorbing 
 time homogeneous Markov 
 chain and $r:S\rightarrow {\bf R}$ a harmonic 
 function.   For every non-absorbing 
 $s\in S$ let  $p_s$ be an alternative transition 
 at $s$ such that $g(p_s)>0$.
 Define a new time homogeneous Markov chain according to the 
 $p_s$. Let   $ a_*: S\rightarrow [0,1]$ be the absorbing rates 
 corresponding to the new time homogeneous Markov chain and let 
 $ r_*:S\rightarrow {\bf R}$ be a  harmonic function with respect 
 to the new transitions such that 
 $r_*$ agrees with $r$ on the absorbing states.  If
 $|v^r(p_s)-r(s)|\leq \delta_s$
 and $ a_*(s)\geq \epsilon_s g(p_s)$ for $0<\epsilon_s\leq 1$ 
  and 
 all non-absorbing $s\in S$ (with $g(p_s)=a(s)$ if  $p_s$ was the  
  original  transition at $s$)
 then the new Markov chain is absorbing and  
 $| r_* (s)-r(s)|\leq \sum_t \delta_t/\epsilon_t$ for all states $s$. 
\vskip.2cm 

{\bf Proof:} The new Markov chain is absorbing because $a_* (s)>0$ 
 for all $s\in S$. With a start at any state $s_0$, we can bound the change 
 $| r_*(s_0)-r(s_0)|$ by the sum over all states $t\in S$ 
 of the one stage deviation at $t$ multiplied by the 
    expected number of visits to the state $t$. 
The  deviation   from one visit to a state $t$ 
 is  bounded by  $| w^r(p_t)-r(t)|$, and since   
 $1/ a_*(t)$ is the expected number of visits to the state $t$ we 
 have the total deviation bounded by $\sum_t
{| w^r(p_t)-r(t)|\over  a_*(t)}$.  
  $| w^r(p_t)-r(t)|\leq g(p_t) |v^r(p_t)-r(t)|$ 
 implies  $| w^r(p_t)-r(t)|/ a_*(t)\leq |v^r(p_t)-r(t)|/\ep _t$. 
\hfill $\Box$\vskip.2cm

\section {Changes in Taboo Probabilities}

In all the lematta of this section, $S$ is a finite state space 
 of a time homogeneous Markov chain. 

\subsection{Reaching a State}

                     For the first three lemmatta we look at 
                      what happens when  a  
 fraction of $P^t(t,s)$ is removed 
                       from the transitions at all  $t$  in a 
 set $T$. 
\vskip.2cm 

{\bf Lemma 3.1}
Let $s$ and $t$ be two  distinct
states and $A$ and $B$ two subsets of states 
 such that $A$, $B$ and $\{ s,t\}$ are mutually disjoint.      
Let $p$ be a part of the   
 transition at $t$ such that  at least positive  $\gamma<1$ 
 of the transition $P^{B\cup \{ t\}}(t,A)$  
   goes through $p$ (meaning that if the complement of $p$ were 
 removed and replaced by motion that went back to $t$ on 
 the next stage with certainty then the new quantity for 
$P^{B\cup \{ t\}}(t,A)$ 
 would be at least $\gamma$ times the old quantity). 
 If the existing transition at $t$ were replaced by $p$ (followed 
 by normalization) and the new transitions were indexed by 
 $*$ then  $P^{B\cup \{ t\}}_*(t,A)\geq \gamma P^{B\cup \{ t\}}(t,A)$ and 
 $P_*^{B\cup \{ s\} }(s,A)\geq \gamma P^{ B\cup \{  s\}}(s,A)$.
\vskip.2cm   
 
{\bf Proof:} $P^{B\cup \{ t\}}_*(t,A)\geq 
\gamma P^{B\cup \{ t\} } (t,A)$ is given. 
If there was never motion from $s$ to $t$ or from $t$ to $s$ then 
 the inequality   $P_*^{B\cup \{ s\} }(s,A)\geq \gamma 
P^{B\cup \{ s\}}(s,A)$ would 
 also be straightforward.  So let us assume that  there is  
  some motion in both directions between $s$ and $t$, and let $A'$ be 
 the set $A$ unioned with all the other states from which there 
 is no motion to either $s$ or $t$. 
 
To estimate $P^{B\cup \{ s\} }_*(s,A)$ let 
    $b:= P^{B\cup \{ s,t\} } (s, A)$, $c:= P^{B\cup A'\cup \{ s\}} (s,t)$, 
     $d:= P^{B\cup A'\cup \{ t\}} (t,s)$
 and $e=P^{B\cup \{ s,t\}} (t,A)$.   Let $   d_*$ 
      and $ e_*$ stand for the contributions to $d$ and $e$ made by  
       the transitions in $p$,
 so that $   d_*\leq d$ and $   e_*\leq e$. 
        By assumption we have  
        $    e_*  +    d_* {b\over b+c}
        \geq \gamma (e+ d{b\over b+c})$. 
         We suppose for the sake of contradiction that 
$\gamma  P^s(s,A)= \gamma (b+{ce\over d+e}) 
> b + {c   e_*\over    d_* +    e_*}=   P_*^s(s,A)$.  
Re-write as $ ( d_*+    e_*) (b   e_* +  b   d_* + c   e_*)> 
(b  e_* +  b   d_* + c   e_*) (d+e)$ or 
 $   d_* +   e_* > d+e$, a contradiction. 
\hfill $\Box$
           \vskip.2cm

{\bf Lemma 3.2:}
    Let $T$ and $A\cup U$ be mutually disjoint  
     subsets of $S$.
      If no more than a frequency of 
       $\gamma  P^{T\cup \{ u\}}(u,A)$ is removed from the transitions of all
        $u\in U\backslash A$
  for some  fraction $0< \gamma <1/(2|U|)$   and no more than 
  a frequency of $\gamma$ in the case of  
   $u\in U\cap A$,                                
        followed by normalization, then 
           for all $x\in S\backslash A$   
           the new resulting probabilities
 $  P_*^{T\cup \{ x\}} (x,A)$ 
               satisfy 
   $P_*^ {T\cup \{ x\}} (x,A)\geq (1- \gamma |U|)
 P^ {T\cup \{ x\}}(x,A)$ and 
     for every $a\in A$  and $x\not\in A\cup T \quad $
   $(1-3|U|\gamma) P_*^{T\cup A}(a,x)\leq P^{T\cup A}(a,x) $.  
              \vskip.2cm

    {\bf Proof:} 
  For $U=\emptyset$ there is nothing to prove. 
       Now assume the result for $U\backslash \{ u\}$, and let   
        $ P_+$ stand for the probabilities where the changes are 
         made in $ U\backslash \{ u\}$.
       Since by induction 
       $  P_+^ {T\cup \{ u\}} (u,A)\geq (1- \gamma |U|+\gamma) 
       P^ {T\cup \{ u\}}(u,A)$, 
          the frequency removal  at $u$ 
         is no more than 
        ${\gamma\over1- \gamma |U|+\gamma}  P_+^ {T\cup \{ u\}} (u,A)$.  
        By Lemma 3.1  
         applied to the case of only one change at $u$, we have 
           for all $x$        
           $P_*^ {T\cup \{ u\}} (x,A)\geq (1- {\gamma\over  
           (1-\gamma |U|+\gamma)}) 
           P_+^ {T\cup \{ x\}}(x,A)\geq (1- {\gamma\over  1-\gamma |U|+\gamma})
            (1-\gamma |U|+\gamma)) P^ {T\cup \{ x\}}(x,A) =
            (1-\gamma |U|)P^ {T\cup \{ x\}}(x,A)$.   

For the second half, if $u\in A$ then it follows by induction because 
 the only way to increase this probability is 
 through the normalization.  Otherwise
        express $P_*^{T\cup A} (a,x)$ as $P_*^{T\cup A\cup \{ u\}}(a,x)+ 
        {P_*^{T\cup A\cup \{ x\}}(a,u)P_*^{T\cup A\cup \{ u\}}(u,x)
        \over 1-P_*^{T\cup A \cup \{ x\}}(u,u)}$.  We notice that 
        $1-P_*^{T\cup A\cup \{x\}}(u,u)
        \geq P_*^u(u, T\cup A\cup \{ x\})\geq P_*^{T\cup \{ u\}} (u,A)$, 
        so that  the change   
        $1-P_+^{T\cup A \cup \{ x\} }(u,u)$  
        to $1-P_*^{T\cup A \cup \{ x\}}(u,u)$ cannot be a decrease by 
 more than a factor of $\gamma / (1-\gamma |U| +\gamma)\leq 2\gamma$.  
                                   The rest follows by 
        $(1-\gamma)P_*^{T\cup A \cup \{ u\}}(u,x)
        \leq P_+^{T\cup A \cup \{ u\}}(u,x)$, (since the only 
 way to increase this probability is through the normalization).  
                                \hfill $\Box$             \vskip.2cm

{\bf Lemma 3.3} Let $T$ be a subset of $S$ and let  
  $s$ be a fixed state such that  $s$ is reached with positive probability 
 from every $t\in T$.  For 
 every $t\in T$ let $q^t$ be a  part of the transition 
  at the state $t$ satisfying  
 $f_{q^t} P_{q^t}^t(t,s) \leq \gamma P^t(t,s)$ where 
 $P_{q^t}^t(t,s)$ is the resulting 
 taboo probability if $q^t$ is a replacement 
 transition at $t$.  Consider  new transitions resulting from 
 the removal of  
 the part $q^t$ at every 
 $t\in T$, followed by normalization.
 If  $|T|\gamma <1$ then 
 $s$ is also reached with positive probability from all of 
 $T$ after the changes. 
\vskip.2cm 

{\bf Proof:}   
We prove  by induction on the size of $T$; by Lemma 3.1 the 
 claim holds for $|T|=1$.  With $v\in T$ also fixed, let us assume 
 that there is some state $u\in T$ such that after the changes
 from a start at $v$ the state 
 $u$ is not reached at all. 
Whether or not 
            one reaches $s$ from $v$ with  the changes  cannot 
            not be influenced by any change made at $u$.  Therefore 
 by the induction hypothesis,  considering changes made in the  smaller set 
$T\backslash \{ u\}$, we have our result. 

Now assume that with the changes all member of $T$ are reached 
 from $v$.   For every pair $t,u\in T$ let $w_t(u)$ be 
 the probability in the original 
 Markov chain with respect to a start at $t$ 
 that $s$ is reached  and that the last visit to a state in $T$ 
 was at the state $u$. Because starting at $t$ rather than at $u$ 
 cannot be a better way to reach $s$ through the state $u$, we have 
 $w_t(u)\leq w_u(u)$.  But then there must 
 be a $u\in T$ such that $w_u(u) \geq {1\over |T|} \sum_t w_t(t)\geq 
{1\over |T|} \sum_t w_u(t)$. This means that at least $1\over |T|$ of the
 original 
 motion $P^u(u,s)$ went directly to $s$ without passing through 
 any other member of $T$ (and therefore after the changes there is still 
motion from $u$ to $s$).  
 \hfill $\Box$ \vskip.2cm

     The following lemma concerns  transitions 
 in two person stochastic games, but  can be generalized   
      to any time homogeneous Markov Chain whose transitions are determined 
 by two independent variables.   
\vskip.2cm

           \vskip.2cm 
          {\bf Lemma 3.4}
          Let $R$ be a subset of non-absorbing states, $U$ a subset of $R$,  
       and $(x,y)$ a pair of stationary  
        strategies such that there is some motion between 
 all pairs of states in $ R$. Let $s,t\in U$ be  special states.          
         Assume for every $u\in U\backslash \{ s\}$  
         that no more  than a frequency of $\gamma P^u(u,s)$ is removed from 
         $x^u \in X^u$  and 
         no more than a frequency of $\gamma$ from  
         $x^s$, followed by normalization; let $\overline x$ stand for  
           the result.
          Assume             for the state $t\in U$    that 
           $P^s_{ \overline x,y} (s,t)\geq \ep P_{x,y}^s(s,t)$. 
 Let $y_*^u$ be a part of $y^u$ for any $u\in U$ 
 with $f_*^u$ its frequency.            
Assume  for all $u\in U$ and both $z\in \{ s,t\}$ that         
            $f^u_* P^u_{x, (y| y^u_*)}(u,z )\leq \delta 
P_{x,y}^u(u,z)$ where  $(y|y^u_*)$  is 
 the strategy that is $y^v$ when $v\not= u$ and is 
 $y^u_*$ otherwise.  
              Let $\overline y$ stand for the result when 
 $y^u_*$ is removed from $y^u$ for every $u\in U$, followed by normalization. 
               Given 
              that $(1-4\gamma |U|) \ep > \delta |U|$   
           with $(\overline x, \overline y)$ there is  
 some motion  from all states in $R$ to $s$ and also some motion  
               from $s$ to $t$.  
               \vskip.2cm

         {\bf Proof:}  
         Since the part of $P_{x,y}^u(u,s)$ that was removed 
 cannot exceed $\delta + \gamma$ of the whole, 
 we have from Lemma 3.3 that  
         $s$ is reached from  all states $v$ in $R$.

As with the proof of Lemma 3.3 we can assume by induction 
 that  all $u\in U \backslash \{ t\}$ are reached 
           from $s$ with 
           $\overline x$ and $\overline y$. 
   We account for $P_{x,y}^s(s,t)$ by considering  
   the last state visited on the way from 
    $s$ to $t$.   
For any choice of $(\tilde x, \tilde y)$ 
 let  $p_{\tilde x,\tilde y}(u,t):=P_{\tilde x,\tilde y}^U(u,t)$ 
be the probability of 
     moving from $u$ to $t$ with no other member of $U$ in between. 
      Let $U':= U\backslash \{ s,t\}$.
     We have  
     $$P_{\tilde   x,\tilde   y}^s(s,t)= p_{\tilde   x,\tilde   y}(s,t)
     +\sum_{u\in U'}{p_{\tilde   x,\tilde   y}(u,t) 
     P_{\tilde   x,\tilde   y}^{\{ t, s\}} (s,u)\over  
     1-P_{\tilde   x,\tilde   y}^{\{ s,t, u\} }
     (u,u)},$$  since  
     ${ P_{ \tilde   x,\tilde   y}^{\{ t, s\}} (s,u)\over  
     1-P_{\tilde   x,\tilde   y}^{\{ s,t, u\} }
  (u,u)}$ is the expected number of times that $u$ is visited before reaching  
      $t$ or returning to $s$,
     with 
     $1-P^{\{ s,t,u\}}(u,u)= P^u ( u,\{s,t\}\cup A)\geq 
P^u(u,s)$, where $A$ is the set from which there is no motion to the set $R$.  
    Define for all  $u\in U'$   
    $e(u):= {P_{ x, y}^{\{ t, s\}} (s,u)\over  
     1-P_{ x,y}^{s,t, u}
     (u,u)}$, with $e(s)=1$, and 
 define $e_*(u)$ correspondingly 
 with respect to $\overline x$ and $y$, with $e_*(s)=1$.
   By Lemma 3.2 we have 
 $(1-4\gamma (|U|-2)) e_*(u)\leq e(u)$ for all $u\in U'$.          
   We can conclude  
    that $$ \sum_{u\in U\backslash \{ t\}} 
   e(u) p_{\overline x,y}(u,t)\geq
    (1-4(|U|-1)\gamma) P^s _{\overline x,y} (s,t)\geq$$ 
   $$\ep (1-4(|U|-1)\gamma) P^s_{x,y}(s,t)= 
   \ep (1-4(|U|-1)\gamma)  
     \sum_{u\in U\backslash \{ t\} } e(u) p_{ x, y}(u,t).\spacenine (4)$$   
Next  define 
$p_{x,\overline y}(u,t):= P^{U} _{ x,\overline y} (u,t).$
 By recognizing that $p_{x,y}(u,t)e(u)/P^s(s,t)$, the probability 
 that the last visit to $U$ was at $u\in U$ from a start 
 at $s$, is less than or equal 
 to the probability that the last visit to  $U$ was 
$u$ {\em with a start at} $u$ (both according to $(x,y)$),
 we have from the defining 
 condition on $\overline y$ that 
$|p_{x,y}(u,t)e(u)-p_{x,\overline y}(u,t)e(u)|\leq \delta P^s(s,t)$.
 After summing over $U\backslash \{ t\}$ 
  we get $$   \sum _{u\in U\backslash \{ t\}}
 e(u) p_{x,\overline y}(u,t)\geq (1-\delta |U|+\delta) \sum_{u\in U\backslash 
 \{ t\}} 
 e(u) p_{x,y}(u,t) \spacenine \quad \quad (5).$$ 
    To show that $u$ reaches $t$ for some $u\in U\backslash \{ t\}$,
 it suffices 
     to show that $p_{\overline x,y}(u,t)+ p_{x,\overline y}(u,t)>p_{x,y}
     (u,t)$ for some $u\in U\backslash \{ t\}.$  But  
      assuming that $p_{\overline x ,y}(u,t)+p_{x,\overline y}(u,t)
      \leq p_{x,y}(u,t)$ for all 
      $u\in U\backslash \{ t \}$, from    
       the above sums in (4) and (5) we must conclude 
        that $1-\delta |U|+\ep (1-4\gamma |U|) 
< 1$, a contradiction to      the initial assumption. \hfill $\Box$
\vskip.2cm 

\subsection{Continuity and Exiting}

Because of the unlimited number of stages, taboo probabilities and harmonic 
 functions of time homogeneous 
 Markov chains are not continuous with respect 
 to  absolute changes in transition probabilities. However, there is 
  a continuity  for relative changes in these transitions.  A  
 result of the same spirit but in a different formal  context 
 is contained in Freidlin and  and Wentzell (1984).   
\vskip.2cm

{\bf Lemma 3.5} 
 Assume that the transitions 
$p^s\in \Delta (S)$ at a subset $U$  
  are changed  
 such that for all $t\in S$, including 
 $s=t$, the resulting $p_*^s(t)$ differs from $p^s(t)$ by 
 no more than a factor of positive $\gamma< 1/(2|U|)$ (necessarily with   
 $p_*^s(t)=0$ if and only if $ p^s(t) =0$).   
Let $ P_*^T(s,A)$ stand for the resulting 
 taboo probability.  For all choices 
   of $s$, $T$, and $A$ with 
 $T\cap A=\emptyset$,
  $P^T(s,A)$ differs from 
 $ P_*^T(s,A)$ by a factor of at most $4 \gamma |U|$. If the 
 original Markov  is absorbing then the resulting Markov chain is 
 absorbing and if 
  $r:S\rightarrow {\bf R}$ is a harmonic 
 function with respect to the original Markov chain  and 
 $ r_*$ is the resulting harmonic function 
 that agrees with $r$ on all the absorbing 
 states then 
   $|r(s)- r_*(s)|\leq 4\gamma |U|M$ for every $s\in S$, 
where  $M$ is a bound on the difference between the function 
 values of $r$ on these absorbing 
 states.    
\vskip.2cm 

{\bf Proof:}  Let $U:= \{ s_1, \dots , s_N\}$.
   Let $P_i^T(s,A)$ stand for 
 the taboo probability when the changes are made only at the subset 
$\{ s_1, s_2, \dots , s_i\}$, and define $\esc_i(t,s)$ in the same way.

First we claim  that for every fixed choice 
 of $s,T,A$ with $s\in T$ that  $P_i^T(s,A)$ and 
 $P_{i-1}^T(s,A)$ differ  at most by a factor of $2\gamma$.    
 Since both $P_i^T(s_i,A)$ and $P_{i-1}^T(s_i,A)$ are expectations 
 over the next stage of some probabilities,     
 we have our claim for $P_i^T(s_i,A)$ and a factor 
 of $\gamma$ by the defining 
 assumption. If $s\not=s_i$ 
then  we get our result from the same observation and 
 the formula 
$P_i^{T} (s,A)=P_i^{T\cup \{ s_i\} } (s,A) + P_i^{T}(s,s_i) 
        P_i^{T\cup \{ s, s_i\} }(s_i,A)
        / P_i^{s_i}(s_i,T\cup B\cup A\cup \{ s\}) $, where 
 $B$ is the set of states such that in either the $i$th or $i+1$st 
 Markov chain there is no motion to the state $s_i$ from the set $B$.

   From formula (1) we have  
 $1- {\esc}_N(t,s)=  P_N ^s(s,t)/( P_N ^s(s,t)+ 
 P_N ^{\{ s,t\} } (s,B))$ and from above  
 that $1- {\esc}_N (t,s)$ does not differ from $1-\esc (t,s)$ by 
 more than a factor of $2\gamma N$.
  Notice that $1-a(s)$ can be written 
 as the expected value of $1-\esc (t,s)$ on the next stage, and 
 therefore  $1-a(s)$ does not differ from $1- a_N(s)$ by more 
 than a factor of $2\gamma N$, where $ a_N$ is the resulting 
 absorption rate.  This implies that   $a(s)=1$ if and only if  
 $ a_N(s) =1$ and in this case we have 
         $ P_N^T(s,A)=  P_N^{T\cup \{ s\} }(s,A)$, 
$ P^T(s,A)= P^{T\cup \{ s\} }(s,A)$, and our result.
Given $a(s)\not= 1$ then by 
 $P^T(s,A)= P^{T\cup \{ s\}}(s,A)/(1- a(s))$  
and $ P_N^T(s,A)=  P_N^{T\cup \{ s\}}(s,A)/ (1- a_N(s))$ 
 we also have our  result.  The claim 
 concerning harmonic functions follows by considering 
 $A$ to be any subset of absorbing states. \hfill $\Box$\vskip.2cm

Next we  define the concept of  exit. (Due to the lack of the semi-algebraic 
 analysis, we will be more restrictive in our definition of an exit 
 than Vieille 2000a or Solan 2000.)  
  For any subset $P$ of non-absorbing states   
a system of {\em exits} from $P$ is a collection of parts of the transitions 
 at the states in $P$ such that all motion from $P$ to $S\backslash P$ 
 must occur through one of these parts. Each part in the collection 
 is called an exit. Given that the Markov chain is absorbing any 
 subset of non-absorbing states must have a system of exits.

  Assume  that there is 
 a partition ${\cal P}$ of the  states such that $\{ s\}$ is 
 in ${\cal P}$ for every absorbing $s$ and  for every non-absorbing 
 $s\in P\in {\cal P}$
   $q^s\in \Delta (S)$  is the transition defined conditionally by the 
union of   
  all the exits from $P$ at the state $s$. Let $A$ be 
 the set of absorbing states.     
  For every $P\in {\cal P}$ let
 $s_P\in P$  be a  representatative for the set $P$. 
   We will create two new time homogeneous 
Markov processes, 
 one  by extending  the state space and the other by 
 contracting it. These constructions are also in Vieille (2000c). 

First we extend the state space. For every $s\in P\in {\cal P}$, create 
 two new states $s^a$ and $s^b$.  Define $S_*:= \{ s^a\ | \ {s\in A}\} 
\cup_{s\in S\backslash A} 
 \{ s^a, s^b\} $, and the corresponding Markov chain 
 will be indexed by $*$.  The  states $\{ s^a\ | \  s\in A \}$ 
 remain absorbing.
 At $s^a$ with $s\in P\in {\cal P}$,
 the motion  goes deterministically 
 to $s_P^b$.  At $s^b$  the transition is 
 labeled  
$p_*^{s^b}\in \Delta ( S_*)$. Let $f_s$ be frequency 
 with which  $q^s$ is used.   
  Let 
 $\overline p^s$ be the transition defined by $p^s$ conditioned 
 on the non-use of $q^s$, given of course that $f_s\not=1$. 
 Define  $p_*^{s^b}(t^a)=  
       f_s q^s (t)$ and 
 $p_*^{s^b}(t^b)= (1-f_s) \overline p^s (t)$ (and otherwise 
    zero if $\overline p^s$ is not defined), with 
 $p_*^{s^b} (a)=p^s(a)$ if $a\in A$.   

Given that the Markov chain  is absorbing, 
next we contract the state space. Define $S_{\sharp}=  \{ s_P\ | \ 
P\in {\cal P}\} $. 
 A previously absorbing state remains absorbing. For every non-absorbing  
 state $s_P$ let the transition at $s_P$ be induced by the distribution 
 on  the next 
 state $t^a$   following 
  $s^b_P$ in the above Markov chain defined on $S_*$. If 
 $t^a$ is absorbing, then $t$ is that next state. If 
 $t^a$ is not absorbing, the $u=u_{P'}$ is the next state 
 with $t\in P'\in {\cal P}$.  Since the  Markov chain on $S_*$ is absorbing, 
 modulo events of zero probability  the transitions of $S_{\sharp}$ 
 are well defined. In a different context (without  taboo probabilities)  
 a similar  statement to the next lemma 
was proven by Vieille (2000c).\vskip.2cm

{\bf Lemma 3.6:} Assume that the Markov chain is absorbing. 
Let $r$ be a harmonic function on $S$ and   
    $M>0$  a uniform bound on all differences 
 in the values of   $r$. 
Let $N$  be the number of the ${\cal P}$ that are not  
 singletons,   and let   
 $0< \delta <  {1\over 2N}$ be given. 
  Assume 
 for every $P\in {\cal P}$ and every distinct pair  $s,t\in P$ 
  that the probability of moving 
 from $t$ to $s$ without passing 
 through any exit of $ P$ is 
 at least  $1-\delta$. The new processes on $S_*$ and 
 $S_{\sharp}$ are absorbing and for any 
 pair of subsets $A$ and $T$ that are unions 
 of members of ${\cal P}$ with $A\cap T= \emptyset$ we have 
 that $P_*^{A_*}(s,T_*)$ differs from $P^A(s,T)$ by no more 
 than a factor of $4N\delta$, where
 $B_*:= \{ s^a, s^b \ | \ s\in B\}$ for all subsets $B$. 
   With  
$ r_*$ representing 
 the new harmonic function on $S_*$ determined by the expected 
value of $r$ on the absorbing states and $r_{\sharp}$ the same 
 for $S_{\sharp}$ we have $ r_* (s^a_R) = r_{\sharp}(s_R)$ for 
 all representative states $s_R$  and  
 $| r_*(s^a)-r(s)|\leq  4MN\delta$ for all $s\in S$. 
    \vskip.2cm

{\bf Proof:} 
Define two new transitions $(\hat p^s \ | \ s\in S)$ and 
 $(\overline p^s\ | \ s\in S\})$ on $ S$. 
$\hat p^s$ is determined by the  
 distribution on  the next state $t^a$ 
 in $ S_*$  from a start at $s^a\in  S_*$. 
  $\overline p^s$ is defined likewise, however from a 
 start at $s^b\in  S_*$. The distribution on the    
 states outside of $P$ with the $\overline p^s$ is the same as 
 with the original transitions $p^s$ on $S$. 
 Because of our assumption 
 concerning the   avoiding of  exits, Lemma 3.5 applies to 
 the difference between $\hat p$ and $\overline p$. 
The claim for the taboo probabilities follows directly from 
 Lemma 3.5, as does also the claim for the harmonic functions. 
  \hfill $\Box$ \vskip.2cm

Lemma 3.6 works because it is based upon the rare use of an exit.  Much more
 problematic is analysing  the consequences of 
the certain use of an exit.  This is the content 
 of Lemma 3.7.    
\vskip.2cm 

{\bf Lemma 3.7} Assume the context of Lemma 3.6 and that $p$ is an 
 exit from $P$ at  
  $t\in P$ with $g(p) >0$. We have \newline 
1)    $|g(p)-g_{\sharp}(p)|\leq 
 4N\delta+ \delta$,\newline 
2) $g(p)$ and $g_{\sharp}(p)$ differ by a factor of no more 
than $4N\delta + {2\delta \over \nu(p)}$, \newline 
3) $\nu (p)$ and $\nu_{\sharp}(p)$ 
differ by a factor 
 of no more than  $4N\delta +2\delta + 
{4N\delta +\delta \over g(p)}$, \newline 
4) $|\nu(p)-\nu_{\sharp}(p)|\leq 8N\delta + 4\delta$, \newline 
5)  $|v^r(p)-v^{r^{\sharp}}(p)|\leq M \min \{ 
   {8N\delta + {\delta \over g(p)}}\ ,  \ 
 {8N\delta +  {2\delta  \over \nu(p)}})\} $.
      \vskip.2cm 

{\bf Proof:} 1)  We define  
 $\hat g$ to be the 
 probability that there is no return to the set $P$ after 
 using the exit $p$ in the original Markov chain.
   From Lemma 3.6 
 we see that $\hat g$ is within a factor of $4N\delta$ of $g_{\sharp}(p)$. 
   From  the avoiding of exits we get 
 that $|\hat g-g(p)|\leq \delta$, which suffices.   

2)  By definition $g(p)\geq \hat g$.
First we show that $\esc(u,t)\leq \delta \hat g /((1-\delta) \nu (p))$ 
 for all $u\in P$. 
Define $w_u$ be the probability that  $p$ will be used  
 before returning to $u$ from a  
 a start at $t$ (with $w_u\leq \delta$ for all $u\in P$). 
Define $\nu_u$ to be the probability that the last visit 
 to $P$ is through the exit $p$ from a start at $u$; 
we have 
 $\nu_u\leq w_u \hat g/(w_u\hat g + \esc (u,t))$, 
 which translates to $\esc (u,t)\leq w_u\hat g/\nu_u$. 
Finally notice that $\nu_u$ doesn't differ from $\nu(p)$ by a factor 
 of more than $\delta$.

Next we compare $g(p)$ with $\hat g$. For every $u\in P$ let $\lambda_u$ 
 be the probability that there is  a return to $P$ from 
 the use of $p$ in the original Markov chain and 
 that $u$ is the first member of $P$ reached. 
Notice that $\sum_u \lambda_u = 1-\hat g$.
We have $g(p)= \hat g + \sum _u \lambda_u \esc (u,t)$.
 This suffices for $(1-2\delta/\nu(p)) g(p)\leq \hat g$.   
Now use Lemma 3.6 for the conclusion.     

3) By definition 
$\nu_*(p)= f_p g_*(p)/a_*(t^b)$ and $\nu(p)=f_p g(p)/a(t)$. 
One way to perceive  $a(t)$ is as the reciprocal of 
 the expected number of visits to $t$ from a start at $t$.  
 With this perspective by Lemma 3.6 and the avoiding of exits 
 we get that $a_*(t^b)$ and $a(t)$ don't differ by a factor 
 of more than $4N\delta +\delta$. This means 
 that if $g_*(p)$ and $g(p)$ don't differ by a factor 
 of more than $\gamma$ then $\nu_*(p)$ and $\nu(p)$ don't differ 
 by more than a factor of $\gamma + 4N\delta +\delta$. Since 
 $\nu_{\sharp}(p)$ is also equal to the probability 
 that the last visit to $P$ starting at 
 $s^b_P$ in the Markov chain $S_*$ went through 
 the exit $p$ we have that 
 $\nu_*(p)$ is within a factor of $\delta$ of $\nu_{\sharp}(p)$ 
 and therefore $\nu_{\sharp}(p)$ and 
 $\nu(p)$ don't differ by a factor of more than 
 $\gamma + 4N\delta + 2\delta$. 
By the same argument as in Part 1 comparing $g_{\sharp}(p)$ 
 with $g(p)$ we get  $|g_*(p)-g(p)|\leq 4\delta N + \delta$ and
 therefore  $g_*(p)$ and $g(p)$ cannot differ by a factor 
 of more than $4\delta N +\delta\over g(p)$   
  and 
 our conclusion.  

4) 
The argument of Part 2 can be repeated with the Markov chain 
 defined on $S_*$ instead of the original on $S$.  The quantity 
 $g(p)$ would be replaced by $g_*(p)$ and $\hat g$ would 
 be replaced by $g_{\sharp}(p)$. 
We have $g_*(p) \geq g_{\sharp}(p)$ and 
$g_*(p)= g_{\sharp}(p) +  (1-g_{\sharp}(p)) \esc_* (s^b,t^b)$.

If 
$g(p)\geq g_*(p)$ we need only 
 $g_*(p)\geq g_{\sharp}(p)$ and the conclusion of Part 2 to get 
  $ g(p) \geq (1-2\delta /\nu(p)-4\delta N) g_*(p)$. Combined 
 with  the arguments from Part 3 we have 
 our goal.   
On the other hand, if $g_*(p)\geq g(p)$ we get our result 
 from  
 repeating Part 2 for $g_*(p)$ and $g_{\sharp}(p)$,
the same arguments of Part 3, plus the claim that   
  $(1-4\delta N -\delta) \esc_*(s_P^b,t^b)\leq \esc(s_P,t)$. 

 $\esc_*(s^b,t^b)$ is no more than 
 $(w+w^2+\dots) h_*$ where 
 $w$ is the probability of reaching an exit of $P$ from 
 $s^b_P$ before returning 
 to $t^b$ and the quantity $ h_*$ is the expected value of $g_{\sharp}$ 
 conditioned  on the use of one of  
 these exits.   
 On the other hand we have that $\esc (s_P,t)$ is at least 
 $w\hat h$ where $\hat h$ is  
  the probability of no return to the set $P$ in the original Markov 
 chain conditioned on the use of one of these exits.   
From Lemma 3.6 we have that $\hat h$ and $h_*$ differ by no more 
 than a factor of $4\delta N$.  That $w\leq \delta$ completes 
 the proof of the  claim.

5) From the proof of Part 1  we had 
 that  $\hat g\geq (1-\delta/g(p))g(p)$ and from  
 Part 2 that $\hat g\geq (1-2\delta /\nu(p)) g(p)$. 
 The rest follows from Lemma 3.6.   
 \hfill $\Box$
\vskip.2cm 

Part 4 of Lemma 3.7 is remarkable because the sum of $\nu$ over 
all transitions in a set $P$ will be $|P|$ rather than something 
 close to one.

\subsection {Polarization}

The process described below, of changing the transitions through a convex 
 combination  
 of two transitions, one giving a higher value and the other giving 
 a lower value of a harmonic function,  
  with  the convex combination yielding the 
 same value, we call {\em polarization}.   
\vskip.2cm

    {\bf Lemma 3.8} 
Let $s$ and $t$ be two  non-absorbing 
states of an absorbing Markov chain.

(i) 
Let $p$ be a part of the   
 transition at $t$ such that
$\nu (p)\geq \ep >0$.

(ii) Let $p$ be a replacement 
transition at $t$ such that $g(p)\geq \ep$.

(iii) 
Let $p$ be a transition at $t$
that is  a convex combination of 
 transitions as described in (i) and (ii). 

    In all three above cases, if we 
    replace the transitions at $t$ by $p$, 
in the case of (i) or (iii) using normalization,
     the  resulting  process is absorbing and the 
    absorption rate of $s$ is at least $\epsilon$ times what is was before 
 the changes were made. 
    \vskip.2cm

    {\bf Proof:}    Let  $b$, $c$,  $d$ and $e$
 stand for the same quantities as in 
  the proof of Lemma 3.1, with $A$ the set of absorbing 
 sets and $B$ the empty set. 

    (i) It follows from Lemma 3.1.

 (ii) Let   
 $ a_*(s)$, 
$d_*$ and $ e_*$  be the corresponding  quantities when 
   $p$ is the transition at $t$.  We assume that  
   $ e_* + d_* {b\over b+c} \geq \epsilon$.  
     Suppose for the sake of contradiction that 
      $\epsilon ( b+ {ce\over d+e}) =\epsilon a(s) >  a_*(s) = 
       b+ {c e_* \over d_* + e_*}$.   Then we have 
        $b e_* + c e_* + b d_* \geq (b+c) \epsilon \geq 
           \epsilon ( b+ {ce\over d+e})> {bd_* + b e_* + 
           c e_*\over d_* +  e_*}$.  This implies 
            $d_* +  e_* >1$, also a contradiction. 

(iii) First we must assume that $b<\ep a(s)$, since otherwise there would be 
 nothing to prove.  Let $a_i$, $d_i$ and $e_i$ for $i=1,2$ stand for 
 the  resulting probabilities 
 from (i) and (ii), respectively, and after normalization 
 in the case of (i). With the convex combinations 
 $\tilde   d:= \lambda d_1 + (1-\lambda) d_2$ and 
$\tilde   e:=\lambda e_1 + (1-\lambda) e_2$ being 
 the new transition quantities, we have that our desired result is equivalent 
 to ${\tilde   e\over \tilde   e +\tilde   d} \geq \ep {e\over e+d} + 
{\ep b-b\over c}$.
 But this follows from (i), (ii), and the fact that  
  ${x_1\over y_1}\geq z$ and  ${x_2\over y_2}\geq z$ implies 
 ${\lambda x_1+(1-\lambda) x_2\over\lambda y_1 + (1-\lambda) y_2} \geq z$ for 
    all non-negative quantities $x_i, y_i, z$ and $0\leq \lambda \leq 1$. 
\hfill $\Box$                                                
\vskip.2cm

{\bf Proposition 3.9} Let $r^1$ and $r^2$ be two harmonic functions, 
 and we assume that the   
 Markov chain is absorbing. Let $N$ be the number of non-absorbing states.   
 Let $ 1$ be a uniform bound on all differences 
 in the values of  $r^1$ and $r^2$. Let $w^1$, $w^2$, $v^1$, and $v^2$ 
 stand for $w^{r^1}$, $w^{r^2}$, $v^{r^1}$, and $v^{r^2}$, respectively.
 Let $1/2> \epsilon>\delta>\gamma >0$, with 
 $\delta< { \ep \epsilon^{3N}\over 2 N(3N)^{N}  }$.
 Let $p^*_s$ be a part of the transition at $s$ such that $ w^2(p^*_s)\leq 
 r^2(s)-\epsilon$ (including the possibility that $p^*_s$ 
 is empty). Assume that if $\nu(p^*_s)\geq \gamma$ then
 there is     an alternative transition $p_s$  
  at $s$ such that $w^2(p_s) \leq r^2(s)-\epsilon$,  
 $|v^1(p_s)-r^1(s)|\leq  \delta$, and there exists another part 
 $q_s$ of the transition at $s$ such that $q_s^d$, the complement of 
 the union of  
 $q_s$ with $p^*_s$,  
 satisfies  $(v^2(q_s^d)- r^2(s)) \nu (q_s^d)\leq N \delta/\ep $. 
 For every subset $T\subseteq \{ s\ | \ \nu (p^*_s) \geq \gamma, \ 
 w^2(q_s)> r^2(s) \}$ 
  define a new time homogeneous Markov chain by the transitions 
 at $t\in T$ defined  
 by $\lambda p_t + (1-\lambda) q_t$ with $\lambda$ satisfying
 $\lambda w^2(p_t)+ 
 (1-\lambda) w ^2(q_t)= r^2(t)$ and furthermore
 for every $v\in S\backslash 
T $ the 
 part $p^*_v$ is discarded, followed by normalization.  
Let the subscript $T$ stand for the quantities  determined by the new 
 transitions with the changes in $T$.  \newline 
{\bf Conclusion:} There is a subset $T\subseteq \{ s\ | \ \nu (p^*_s)
 \geq \gamma, \  w^2(q_s)> r^2(s)\}$ such that the
 new process is absorbing and 
 for both $i=1,2$ and 
 all $s\in S$  
$| r_T^i(s)-r^i(s)|\leq \epsilon$
\vskip.2cm 

{\bf Proof:}
 First we consider what happens when the changes are made only at 
 a set $T$ (meaning that 
 the part $p^*_s$ is kept in for $v\not\in T$),
 which we will label with $T,*$. Because $r^2$ 
 remains a harmonic function after the changes are made  
  and there is always a positive probability at all states in $T$  
 that the harmonic function drops by $ \ep$,
   the resulting 
 time homogeneous Markov chain is absorbing with   
 $ r_{T,*}^2(s)=r^2(s)$ for every $s\in S$.

Next we must determine which subset $T$ will be chosen. 
Choose any $t_1$  such that $\nu (p^*_{t_1})\geq \epsilon ^2/2N$,
 and put $t_1$ 
 in $T$. If there 
 exists no such $t\in S$ then let $T$ be the empty set. At any set $T$ with 
 $|T|=k-1$ 
 formed so far,    
 put into $T$ any $t_{k}$ such that $\nu_{ T,*} (p^*_{t_k})\geq 
 \epsilon ^2 / 2N$, and stop if there is no such new state $t_k$.

{\bf Claim:} For  any set $T$ that has been 
 already chosen and any $t\not\in T $ that 
 could be added to $T$ we have 
 $a_{T\cup \{ t\}, *} (u) \geq  { \ep^3 \over 3N} a_{T,*}(u)\geq 
{\ep^3\over 3N}{ \ep ^{3|T|}\over (3N)^{|T|}} a(u)$ for all $u\in S$, 
 $g_{T,*}(q_t)\geq {\ep^3\over 3N}{\ep ^{3|T|}\over (3N)^{|T|}} 
 g(q_t)$ and $w^2(q_t)> r^2(t)$. 

{\bf Proof of Claim:} Assume  
 that $t$ will be added to $T$.
Look at the transition $q_t^d$ and the indentities  
 $w_{T,*}^2(q_t^d)-r_{T,*}^2(t) =w^2(q_t^d)-r^2(t) 
=(v^2_{T,*}(q_t^d)-r_{T,*}^2(t)) g_{T,*}(q_t^d) =(v^2(q_t^d)-r^2(t)) g(q_t^d)$ 
 from the fact that $r^2$ remains 
 the harmonic function. Consider  
 the definitions 
  $\nu_{T,*}(q_t^d)= f_{q_t^d} g_{T,*}(q_t^d)/a_{T,*}(t)$ and 
$\nu(q_t^d)= f_{q_t^d} g(q_t^d)/a(t)$; 
 they show that the new absorption rate determines alone the 
 new value $(v^2_{T,*}(q_t^d)-r^2_{T,*}(t)) \nu_{T,*}(q_t^d)$. 
  From the induction 
 assumption we must conclude that 
 $\nu_{T,*}(q_t^d) (v^2_{T,*}(q_t^d)-r^2(t))\leq 
  { (3N)^{|T|}\over  \ep ^{3|T|}}\nu(q_t^d) (v^2(q_t^d)-r^2(t))\leq 
      { (3N)^{|T|}\over  \ep ^{3|T|}}{N\delta\over \ep} < \ep^3/6N $.
If $q_t^c$ is the union of $q_t^d$ with $p^*_t$ from  
 $\nu_{T,*}(p^*_t) \geq \ep ^2 /2N$ and $w^2_{T,*} (p^*_t)\leq 
 r^2(t)-\ep$ we get 
that $\nu_{T,*}(q^c_t) (v^2_{T,*}(q^c_t)-r^2(t))\leq - {\ep^3\over 3N}$, 
 which implies that $w^2(q_t)> r^2(t)$ and   
 $\nu_{T,*}(q_t)\geq 
  \ep^3/(3N)$.

 Next suppose for the sake of contradiction that 
$g_{T,*}(q_t)< {\ep^3\over 3N}{\ep ^{3|T|}\over (3N)^{|T|}} 
 g(q_t)$.  Since  $\nu(q_t)=
 f_{q_t} g (q_t)/a(t)$,
 $\nu_{T,*}(q_t)= f_{q_t} g_{T,*} (q_t)/a_{T,*}(t)$ and 
 $\nu_{T,*}(q_t)\geq \ep ^3/(3N)$, by the induction 
 assumption  we would be forced to accept $\nu(q_t) >1$, 
an impossibility.

 By Lemma 3.8 we have our claim on the absorbing 
 rates for all states other than $t$. For 
 the state $t$ we have
   $g_{T,*}(q_t) \geq f_{q_t} g_{T,*}(q_t) = \nu_{T,*} (q_t) a_{T,*}(t)\geq 
\ep^3 a_{T,*}(t) / (3N)$.  With     
 $g_{T,*}(p_t)\geq \ep$ our claim is proven.    
\vskip.2cm

With the claim we  conclude from Lemma 2.1 that 
  $|r^1_{T,*}(s)-r^1(s)|\leq {(3N)^{|T|}    \over   
\epsilon^{3|T|}}\delta N
\leq \epsilon/2$ for all $s\in S$.

Next, we must show that it is impossible for any state $s$ not polarized 
 to satisfy   
 $\nu_{T,*}(p^*_s)\geq \epsilon^2/N$.  This holds 
 for all states  with $\nu (p^*_s)\geq \gamma$, by construction. 
Let's assume that   $\nu (p^*_s)< \gamma$;  this means that 
 the probability of ever using $p^*_s$ in the 
 original Markov chain cannot   exceed $\gamma /\ep$. 
But by the above claim we know additionally that the probability 
 of using $p^*_s$ in the altered Markov chain indexed 
 by $T,*$ cannot exceed ${\gamma \over \ep}\ {\ep^{3(N-1)}\over 
 (3N)^{N-1}}< \ep^2/2N$. 

Next we must consider the influence of the removed  
 $p^*_t$ in the above Markov chain indexed by $T,*$.  
  For any $s$ with $\nu_{T,*}(p^*_s)\leq \epsilon^2/2N$   
  the chance of ever using the transition $p^*_t$ cannot 
  exceed $\epsilon /2 N$, and so 
 they cannot contribute an average of more than $\epsilon/2$ to 
either the function $r^1$ or $r^2$. 
\hfill $\Box$

\section {From Markov Chains to Equilibria}
\subsection {Application of the Doob-Kolmogorov Inequality}
We must prove Proposition 4.2, a cornerstone of our analysis.
\vskip.2cm

{\bf Lemma 4.1:}
 Let $X$ be the finite state space of a time homogeneous Markov chain with 
 probability transitions 
 $(p^x\in \Delta (X)\ | \ x\in X)$.
 Let   $v:X \rightarrow 
 {\bf R}$    be a harmonic function and let 
 $M>0$ be a  
 bound for the  maximal difference between all values of  $v$.

For every $x\in X$ define the non-negative quantities 
 $w(x)$   by $w(x) =
\sum_{ y \in X}        p^{x} (y)   |v(y)-v( x)|$. Let $n$ be 
 the number of states $x$ such that $w(x)>0$. 
For any path $p=(x_0, x_1, x_2, ...)$ in $X$ define $w(p)=
\sum_ {i=0}^{\infty}  w (x_i)$.  

{\bf Conclusion}: The expected value of the function 
 $w$ does not exceed $Mn$.   \vskip.2cm

{\bf Proof:} We isolate the problem, handling each state $x$ separately. 
Since $|v(y)-v(x)|$ is always less than or equal to $M$ times
 $\esc (y,x)$,  
 we have that $w(x)\leq a(x)M$. Therefore the  
 part of the sum that comes from visits to $x$ does not exceed 
 $a(x)M\sum _{i=0}^{\infty} (1-a(x))^i =M$.  \hfill $\Box$ \vskip.2cm

{\bf Proof of Proposition 4.2} (as stated in the introduction):
 
Define the random variable $r_i$ on the odd steps $i$ to be 
$v(y_{i})-v(x_{i-1})$, and $R_i$ to be the sum of 
 the $r_k$ for odd $k\leq i$. For $y\in Y_x$ define 
 $r(y)$ to be $v(y)-v(x)$.   

Define  a new quantity, $\tilde   w(x):= \sum_{y\in Y_x} p^x(y) 
|v(y)- v(x)|$. Let $w(x)$ be the old 
 quantity  on the Markov chain from Lemma 4.1 
 defined only on the $X$, -- we  ignore the visits to the $Y_x$ sets, 
 and consider only the motions from $X$ to $X$.   

The Doob  submartingale inequality states that if $(S_i\ | \ 
i=0,1,\dots, n)$ is a martingale with zero expectation then for every 
 $n\geq 0$, positive value $c>0$  and  exponent $p\geq 1$ 
 the probability that $\max_{i\leq n}|S_i|>c$ is 
 less than  ${\bf E}(|S_n|^p)/c^p$ (Williams 1991, Section 14.6).
  Since the martingale 
 property implies that ${\bf E}(S_n^2)$ is equal to the sum over 
all the stages $1\leq i\leq n$ 
of $E(s_i^2)$ where $s_i$ is the change in value between the $i-1$st 
 stage and the $i$th stage, we have for every finite even and positive $Q$ 
$$\mbox {Probability } \Big(\max_{ i< Q } | R_i| > 
  \epsilon \Big)
 < {1\over \epsilon^2} {\bf E}\Big( \sum_{ i< Q,\ y\in Y_{ x_{i-1}}} 
 p^{ x_{i-1}} (y)r(y)^2\Big). $$ 
By taking the limit as $Q$ goes to infinity and $\delta \leq |r(y)|$ we get 
$$\mbox {Probability } \Big(\max_{ i< \infty } | R_i| > 
  \epsilon \Big)
 < {1\over \epsilon^2} {\bf E}\Big( \sum_{ i< \infty ,\ y\in Y_{ x_{i-1}}} 
 p^{ x_{i-1}} (y)r(y)^2\Big ) \leq $$ 
 $$\delta {1\over \epsilon^2}
  {\bf E} \Big ( \sum_{i< \infty,\ y\in Y_{ x_{i-1}}} 
 p^{ x_{i-1}} (y)|r(y)|\Big )=  \delta  {1\over \epsilon^2}{\bf E}\Big(
  \sum_{i< \infty ,\ y\in Y_{ x_{i-1}}}
 \tilde   w(x_{i-1})\Big)  .$$ 
 Since  
 by the triangle inequality $ \tilde   w(x)\leq w(x)$ for all $x$, we have 
$$\mbox {Probability } \Big(\max_{ i< \infty } | R_i| > 
  \epsilon \Big)
 < \delta {1\over \epsilon^2} {\bf E}\Big( \sum_{i<
 \infty,\ y\in Y_{ x_{i-1}}} 
  w (x_{i-1})\Big),$$
and by Lemma 4.1 this is no more than $\delta Mn/\ep^2$. 
  So with $\epsilon \leq 1/2$, 
 we have our result from the size of $\delta$. 
 \hfill $\Box$ \vskip.2cm 

The problem of extending 
 Proposition 4.2 to Markov chains that are not time homogeneous (or have 
 countably many states) 
 lies with    Lemma 4.1 and not in the proof of Proposition 4.2.  

The following corollary relates the above work on Markov chains 
 to    our two-person stochastic games. 
   Because 
 the application of this corollary involves an 
 altered state space,  this result 
 should be understood in an abstract way.\vskip.2cm 

{\bf Corollary 4.3:} 
 Let  $(x,y)\in X \times Y$ be 
  stationary  absorbing  strategies.    
Assume that \newline 1) for both players $k=1,2$ and  $s\in S$ $r_{x,y}^k(s)$
 is greater than  $j_z^k(s)-\epsilon$ with $z=x$ if $k=2$ and 
 $z=y$ if $k=1$, 
  and that\newline   
 2) for both player $k=1,2$ and  all moves $c$ used with positive 
 probability with  
$(x,y)$ by Player $k$ the value $w^{k}_{x,y}(c)$ is within $\delta$ 
of $r_{x,y}^k(s)$. \newline  
 {\bf Conclusion:} For any positive $\epsilon < 1/2$
if $\delta$ is no more than $\epsilon^3\over n$ then    
the strategies  $(x,y)$  generate a 4$\epsilon$-equilibrium of the 
 stochastic game. \vskip.2cm 

{\bf Proof:}
 We define the following strategy for Player $k$.
   For every starting point   
  $ s_0\in   {\cal S}$  
let $n_{s_0}$ be large enough   such that with a start at 
 $s_0$ and the play according to  $(x,y)$ 
 the  probability  that there is no absorption before
 the  $n_{s_0}$th stage is less than $\epsilon /10$. 
Let   $s_0, s_1, \dots $ be any sequence 
 of states reached in the game and for both $k$
 let $c^k_0$, $c^k_1 , \dots$ be the sequence 
 of moves made by Player $k$. For   $k'\not=k$  as long as  
  $\sum_{i=0}^{l} 
(w^{k'}_{x,y}(c^{k'}_i)-r_{x,y}^{k'}(s_i))\leq \epsilon$ and  
 the stage $l$ does not exceed $n_{s_0}$  and 
 Player $k'$ never chooses  $c^{k'}_i$ outside of  the support set 
 of his stationary strategy,   
 then 
 Player $k$  continues to act according to his stationary strategy.
  As soon as one of the  
 above conditions is violated at some stage $l$ then on the next stage 
 $l+1$  
  both players punish eachother according 
 to the functions $c_1 +\ep$ and $c_2+\ep$.  (The mutual punishment is 
 necessary because otherwise a player could intentionally prolong 
 the game with an interest in punishing the other player.   
   The result can be extended to 
 multi-player stochastic games if it can be determined who should punish 
whom in all situations!) 
 That no player $k$  can obtain an expected payoff more than 
 $2\epsilon$ above the function $r^k$ by choosing a different 
strategy    is self explanatory. 
 That punishment occurs before absorption with probability no more than 
$2 \epsilon$ if both players adhere to the suggested strategies follows 
 from Proposition 4.2.
  \hfill $\Box$ \vskip.2cm 
  
\subsection{Situations}

Next we create an expanded state space from the original state space 
 through partitions of the histories. 
For every $s\in {\cal S}$ let ${\cal P}^s$ be a partition of ${\cal H}^s$. 
 Define $\hat {\cal S}$ to be the disjoint
 union  $\cup_{s\in {\cal S}} {\cal P}^s$. 
 For every $t\in \hat {\cal S}$ 
let $b(t)\in {\cal S}$ be the member of ${\cal S}$ such that 
 $t\in {\cal P}^{b(t)}$. A member of $\hat {\cal S}$ we call a {\em situation}.
We define the situations $\hat {\cal S}$ to be {\em normal} if and only 
 if  the next $u\in \hat {\cal S}$ following 
 a  $t\in \hat {\cal S}$ is 
 determined  by the 
situation  $t$, the choice of moves by the players at $t$, and 
  the next $s\in {\cal S}$ with $b(u)=s$.      Normalcy implies that 
 one can define a stochastic game on the situations as a new state space. 
 \vskip.2cm 

 {\bf Corollary 4.4:}     
 Let the situations $\hat {\cal S}$ be normal, let absorbing stationary 
 strategies $(x,y)\in \prod_{s\in \hat {\cal S}} \Delta (A_1^{b(s)})\times 
\prod_{s\in \hat {\cal S}} \Delta (A_2^{b(s)})$ 
 be  defined on the situations $\hat {\cal S}$,  
 with $\hat r_{x,y}^k:\hat {\cal S}\rightarrow {\bf R}$ the expected
 payoff for Player 
 $k$ as determined by the above stationary strategies and the functions  $r^k$ 
 on the absorbing states and $\hat w^k_{x,y}$ the corresponding 
 expected value of $\hat r^k_{x,y}$ on the next stage.   
Assume that \newline 
1) for every $s\in \hat {\cal S}\quad \hat 
r^k(s)\geq j_z^k (b(s))-\epsilon$  where $z=x$ 
 if $k=2$ and $z=y$ if $k=1$ and \newline 
2) for every move $c$ used with positive probability 
 at a situation $s$ by Player $k$  $|\hat w_{x,y}^k(c)-
\hat r_{x,y}^k(s)|\leq \delta$. \newline 
   If $\delta$ is no more 
than $\epsilon^3\over |\hat {\cal S}| $ then these stationary
  strategies  generate a $4\epsilon$-equilibrium of 
 the original stochastic game. \vskip.2cm

{\bf Proof:} Because a stochastic game is defined 
 by the normality of $\hat {\cal S}$ and the conditions of Corollary 4.3 are 
 preserved, the result follows 
 by Corollary 4.3. \hfill $\Box $ \vskip.2cm

\subsection{First Main Theorem}

For any subset $R\subseteq {\cal N}$ and a state $s\in R$, a pair 
 $a\in A^s_1$ and $b\in A^s_2$ of moves is called a {\em primitive exit} 
 from the set $R$
 if with positive probability there is motion from $s$ to $S\backslash R$ 
 using the pair $a$ and $b$.
 By the definition of $\rho$,
 any use of a primitive exit at $s$ results in a probability of 
at least  $\rho$  of reaching the complement of $R$.

  For every subset $B$ of Player Two moves 
 in a set $R$ we define  a {\em $B$
 exit} (or simply exit if there is no ambiguity) 
 from $R$ to be any pair $(a,b)$ of moves at an $s\in R$ such  that 
 $(a,b)$ is already  
 a primitive exit from $R$ or  
 $b\in B$.  Let  $E^B(R)$ stand for 
 the set of all  $B$ exits from $R$. 

Define  $B_{x,y}^{\gamma}(s)$ to be those moves of Player Two at the state 
 $s$ with $w^2_{x,y}(b)\leq  r_{x,y}(s)-\gamma$, and let 
 $B_{x,y}^{\gamma}(R)$ be the union of all the 
 $B_{x,y}^{\gamma}(s)$ for all $s\in R$. 
 For every $s\in {\cal N}$ define  $z^{\gamma}_{x,y}(s)$ to be 
   $\sum_{b\in B_{x,y}^{\gamma}(s)} \nu^b$.
  For any 
 subset $R\subseteq {\cal N}$ define 
 $z_{x,y}^{\gamma}(R) := \sum_{s\in R} z_{x,y}^{\gamma}(s)$.

For any stationary strategy $x\in X $ (or $y\in Y$) 
define a {\em simplication} of $x$ 
 to be another stationary strategy 
$\overline x\in X$ obtained from $x$ by dropping the use 
 of certain moves, followed by   
  normalizing   what remains. 
Call the simplification a $\gamma$-simplication if 
 the frequency  of the moves removed did not exceed $\gamma$. 
 The simplication is {\em within} a set $T$ of states if changes were 
 made only within the set $T$. 
\vskip.2cm

{\bf Theorem 1:} Assume  for every choice of positive 
 $1/2> \ep > \overline \epsilon>  \hat \epsilon > \tilde \epsilon>0$ with 
 $\overline \ep < \epsilon^3/(50|{\cal N}|)$,    
 $\hat \epsilon < {\overline \ep \ \overline \epsilon^{3|{\cal N}|}
\over 5(3|{\cal N}|)^{|{\cal N}|} |{\cal N}|}$ 
 and $\tilde \ep < \overline \ep \ \hat\ep / 40|{\cal N}|$
 that\newline   
{\bf 1)}  there are  absorbing stationary strategies $(x,y)\in X\times Y$ 
 with\newline  
a) $r^2_{x,y}(s)\geq j^2_x(s) -\epsilon /2$ for all
 $s\in {\cal N}$,\newline b)  
 $r^1_{x,y}(s) \geq j^1_y(s)-\ep /2$ for all $s\in {\cal N}$, and\newline 
c) for every move $a$
 of Player One used in $x$ with positive probability  at $s$  
  we have  $|w_{x,y}^1(a)-r_{x,y}^1(s)|\leq 
 \tilde \ep$,  \newline      
{\bf 2)} a partition ${\cal R}$  of  a subset 
 $P\subseteq {\cal N}$   and for every $R\in {\cal R}$ a set $B_R$ of Player 
 Two moves in $R$ containing $B^{\overline \ep}_{x,y}(R)$ 
 such that    
\newline a) $\forall s\not\in P$ $z^{2.5\overline \ep}_{x,y}(s)< 
 \tilde \ep$ and  \newline  
b) for every distinct  $s,t\in R\in {\cal R}$    
  the probability of reaching 
 $s$ from $t$ before using a member of $E^{B_R}(R)$  
 is  at least $1-\gamma^*$ with 
  $\gamma^*:=  \tilde \ep \ \overline
 \ep /(40n|{\cal N}|) $,\newline\newline   
  and for any $R\in {\cal R}$ if $z^{2.5\overline \epsilon} _{x,y}(R)
\geq   \tilde \epsilon $  
 then   there is a  
  special subset $D_R\subseteq R$, 
  a representative $s_R\in D_R$ and  \newline  
{\bf 3)} an $\tilde \ep$ 
  simplication $y_R$ of $y$ within $R$ created by removing 
  the set $B_R$ of moves  such that\newline  
 a)$v^2_{x, y}(b)\leq r^2(s)$ for every $b\in B_R \cap A^2_s$  and  \newline   
b)   $|r^1_{x,y_R}(s_R) -r_{x,y}^1(s_R)|\leq \hat \ep$,\newline   
{\bf 4)}  $\tilde \ep$-simplications 
$(x_{C}, y_{C})$ of $(x,y)$ within  
 $D_R$ such that  with $(x_{C}, y_{C})$ the 
 play never leaves the set $D_R$ and from any state in $D_R$  all other states 
 in  $D_R$ are reached with probability one, and \newline   
{\bf 5)} a strategy $y_D$ created from $y_C$ by adding to $y_C$ in 
 the set $D_R$ small probabilities 
 of using a subset of Player Two moves $V_R$ used in $D_R$ with 
 $V_R \subseteq B^{2.4\overline \ep}_{x,y}(R)$ and 
a real  positive  value $\xi_R\leq r^2_{x,y}(s_R)-2.4\overline \ep$ 
 such that \newline 
a)    with $(x_C,y_D)$ for every pair $s,t\in D_R$ 
 the probability of reaching 
 $s$ from $t$ before using a member of $V_R$  
 is  at least    $1-\gamma^* $\newline    
 b)  $\xi_R \geq j^2_x(t)-\ep$ for 
 all $t\in D_R$, \newline  
c) for all moves $b\in V_R$      
 $|w^2_{x, y}(b)- \xi_R |\leq \tilde \ep $, and \newline 
d)   $ |r^1_{x_C,y_D}(s)-  r^1_{x,y}(s)|\leq 
  \hat \ep $.\newline 
{\bf Conclusion:}  
 With the assumption 
 that Player One can send transition independent 
 signals, 
the stochastic game has approximate equilibria.   
\vskip.2cm 

{\bf Proof:} 
  We define the set $B$ of Player Two moves to 
 be $\cup_{R\in {\cal R}} B_R \cup _{s\not\in P}B^s_2$, and 
 define the exits to  be the $B$ exits. 
Let the corresponding state spaces  ${\cal S}_*$ and 
 ${\cal S}_{\sharp}$ from Lemma 3.6 be induced by $(x,y)$ and 
 the partition ${\cal R} \cup \{ \{ s\} \ | \ s\not\in P\}$.
  For every $s_R\in {\cal S}_{\sharp}$ let 
 $p_R^*$ be the transition at $s_R$ in ${\cal S}_{\sharp}$
 induced by the Player Two moves in $B_{x,y}^{ 2.5 \overline \ep}(R)$. 
 For every $R\in {\cal R}$ 
define $p_R$ to be the alternative 
 transition from $s_R$ in ${\cal S}_{\sharp}$
 induced  by the Player Two moves $V_R$ according  
 to   $(x_C, y_D)$. Define 
  $q_R^c$  to be the transition 
   induced by the moves 
 in $B_R$, and define $q_R^d$ so that 
 $q_R^c$ is  the disjoint union of $q_R^d$ with $p_R^*$.

We will confirm the conditions of Proposition 3.9 on the state space 
 ${\cal S}_{\sharp}$, with $2.4\overline \ep$, $2\hat\ep$, and 
 $2\tilde \ep$ the quantities $\ep$, $\delta$, and $\gamma$ of that 
 lemma, respectively.  

First, by Lemma 3.6 the Markov chain on ${\cal S}_{\sharp}$ is 
 absorbing. 
For $i=1,2$ let $ r_{\sharp}^i : {\cal S}_{\sharp}
 \rightarrow {\bf R}$ be the harmonic
 function  that agrees with   
 the function $r^i$ on the absorbing states. If $\nu_{\sharp} (p^*_{s_R})\geq 
 2\tilde \ep$ then by Lemma 3.7 
$z_{x,y}^{2.5\overline \ep}(R) \geq 3\tilde \ep /2$ 
 and if $s\not\in P$ then $\nu_{\sharp}(p^*_s)\leq 1.1 
 z^{2.5\overline \ep} _{x,y}(s) \leq 1.1 \tilde \ep$. 
 By Lemma 3.6 
 we have for every representative 
 $s_R$ that $ r_{\sharp}^i(s_R)$ is within $4\gamma^* |{\cal N}|$ of 
 $r_{x,y}^i(s_R)$. 
Equally important,  Lemma 3.7 implies that 
 $w^{ r_{\sharp}^2}(p_R^*)\leq  r_{\sharp}^2(s_R)-2.4\overline \ep$,  
  and 
   $ |v^{r_{\sharp}^1}(p_R)-r_{\sharp}^1(s_R)|\leq 11\hat \ep/10$.
  Since $q^d_R$ is induced by some $B_R$ moves  
 by  Lemma 3.7  and Condition 3a we have  
 $(v^{ r_{\sharp}^2} (q^d_R)-   r^2_{\sharp}(s_R)) \nu_{\sharp}(q^d_R) 
< 2\tilde \ep $. 

Left to confirm is that 
$|v^{r_{\sharp}^1}(q_R)-r_{\sharp}^1(s_R)|\leq 2\hat \ep$. 
 We apply Lemma 3.6 to the pair $(x, y_R)$ and the transitions it 
 induces on ${\cal S}_{\sharp}$.  Since the avoiding 
 of exits by  $(x,y)$
 implies the 
 same for the pair $(x,y_R)$, we have that 
 $ | \underline 
r^1 (s_R)-r^1_{x,y_R}(s_R)|\leq   4\gamma^* |{\cal N}|$, 
 where $\underline  r^1$ is the harmonic function induced by 
 $(x,y_R)$ on ${\cal S}_{\sharp}$.  $\underline  r^1 (s_R)$ is equal 
 to $v^{r_{\sharp}^1}(q_R)$.   
 With the given $|r^1_{x,y_R}(s_R)-r^1_{x,y} (s_R)| \leq \hat \ep$ 
 and the above relation of $r^1_{x,y}$ to $r^1_{\sharp}$ we are done 
   establishing the conditions of Proposition 3.9. 

We apply Proposition 3.9 to  $S_{\sharp}$
 with ${\cal T}$ the subset of ${\cal R}$ that has been polarized. 
  We  conclude that 
the new harmonic functions $ \tilder_{\cal T}^i:= (r_{\sharp}^i)_{\cal T}$
 on ${\cal S}_{\sharp}$ 
  satisfy $| \tilder_{\cal T}^i(s)-r^i_{x,y}(s)|\leq 
3 \overline \ep$ for all $s\not\in P$ and 
$| \tilder^i_{\cal T}(s_R)-r^i_{x,y}(s_R)|\leq 
 3\overline \ep$ for all $R\in {\cal R}$.  

Next we define the situations $\hat {\cal S}$,
 with one, two, or three situations 
 defined for each original state in $\cal S$.  For any $s\not\in P$ or 
 for $s\in R\in {\cal R}$ with $R\not\in {\cal T}$ not polarized there is 
 only the situation  $s^e$ (including the case of absorbing states). 
 We always start the game at an $s^e$. 
 At any  situation $s^e$  for $s\not\in P$ or $s$ in  
 a non-polarized $R\not\in {\cal T}$ the players 
 perform $(x^s,\hat y^s)$ where $\hat y^s$ is the $\gamma^*$ simplication 
 of $y^s$ resulting from the removal of all Player Two moves in 
 $B_{x,y}^{2.5\overline \ep} (s)$.
 If  $s$ is in a polarized   
  $R\in {\cal T}$ and is not the representative $s_R $
 the players perform $(x,y_R)$.
 Following any $s^e$ other than 
  $s = s_R$ the next situation is a $t^e$, where 
 $t$ is the next state in $ {\cal S}$. 
Also  following the performance of an exit,
 no matter what the situation 
was on the previous stage, if $t\in {\cal S}$
 occurs on the following stage then 
 the next situation is also $t^e$. 
This means that only motion inside of an $R\in 
{\cal T}$  involves situations other than those with 
 the subscript $e$.

At any  $s\in R\in {\cal T}$  there is 
 either two situations $s^e$ and $s^f$ if $s\not\in D_R$ or 
 three situations $s^e$, $s^f$, and $s^g$ if $s\in D_R$.  
 For such an  $R\in {\cal T}$ 
let $\lambda _R$ be the quantity determined by 
 the  application of Proposition 3.9 to the transitions on 
${\cal S}_{\sharp}$. 
  Since Player One can 
 send signals, for every $s_R\in D_R$ for a polarized $R\in {\cal T}$ 
 we associate one of every pair of her moves with the 
 symbol $f$ and the other with the symbol $g$.  
 If  $s_R^e$ is 
 the present situation  
then with  probability 
 $\lambda_R$ Player One chooses a move associated with 
 the symbol $g$ and with 
 $1-\lambda_R$ a move associated with the symbol $f$; in both cases  
 the players perform $(x_C, y_C)$.
 (Because all moves are paired, we can modify 
  $x_C$ to use only those moves corresonding to $f$ or only moves 
 corresponding to $g$  without changing the transition probabilities in the 
 space $\cal S$.) If $t$ is the next state 
 and a move corresponding to $f$ was used, then $t^f$ is 
 the next situation; otherwise the next situation is $t^g$. 
   At any $s^f$ with $s\in R\in {\cal T}$ the play continues according 
 to $(x,y_R)$, always to a  next  situation  $t^f$ if there was no 
 use of an exit.  On the other 
 hand, from any $s^g$ with $s\in D_R$ the motion follows $(x_C, y_D)$, and 
 unless a move from $V_R$ is used the next situation is a $t^g$, 
 necessarily with $t\in D_R$.

Define $\hat r^i$ to be the harmonic function on $\hat {\cal S}$ determined 
 by the above defined stationary behavior and 
$\hat r^i =r^i$ on the absorbing states.
Given the above conditions, to
 apply Corollary 4.4 it suffices that  neither player $i$ 
 can change the expected value of $\hat r^i$ by more than $10\overline \ep$ 
 at any one stage. 
 With the role of the $\xi_R$ 
  we need only 
 show that $\hat r^i$ is within $\overline \ep$ of $r^i_{\sharp}$ on 
 all the $s_R$ and the $s\not\in P$.  To do this, we introduce two new   
  transitions 
defined on ${\cal S}$, indexed by $\circ$ and ${\cal o}$. $p_{{\cal o}}$ and 
 $p_\circ$ are identical on states $s$ that are not in a polarized $R$, 
 and then it is that induced by the behavior at   the 
 situation $s^e$. At $s$ in a polarized $R\in {\cal T}$ 
  $p_{{\cal o}}^s$ is the distribution  determined 
 by the next situation  $t^e$ following the situation $s^e$. 
  $p_\circ^s$ is determined by the next situation $t^e$ conditioned 
 on having reached either $s_R^f$ or $s_R^g$ before any exit 
 was performed. The $p_\circ^s$ transitions generate harmonic 
 functions $r^i_\circ$ that are identical to $r^i_{\sharp}$ on the 
 ${\cal S}_{\sharp}$, and the $p_{{\cal o}}$ transitions 
 generate harmonic functions $r^i_{{\cal o}}$ that are identical 
 to $\hat r^i$ on the subset $\{ s^e \ | \ s\in {\cal S}\}$.     
 Because $\lambda_R$ cannot 
 be greater than $1-2\overline \ep$ and  the probability from 
 a situation $s^e$  that an exit from 
 the stationary strategies $(x,y_R)$  
 is used before getting to $s^e_R$ is no more than  
 than $\gamma^*$, for every $s,t\in {\cal S}$ the 
transition probability $p^s_{{\cal o}}(t)$ does not differ 
 by more than a factor of $\gamma^*/\overline \ep $ from  
 $p^s_\circ(t)$.   Finally 
Lemma 3.5 implies that the functions $r^i_\circ$ and $r^i_{{\cal o}}$ 
do not differ by more than $4\gamma^*|{\cal N}| /\overline \ep< 
\tilde \ep$.       
\hfill $\Box$ 
\vskip.2cm

   \section {The auxiliary game}

The main issue is to define the
 ``correct'' discounted evaluation of Player Two,
since, as shown in Solan (2000),
 a naive definition of his discounted evaluation 
does not prove equilibrium existence when there are a multitude of 
 non-absorbing states.
    
We assume that positive  $\ep$ and 
 $\overline \epsilon$   have been fixed.

\subsection {The function $\xi$}

Let $b$ be any move of Player Two at a state $s\in {\cal N}$.

For any $(x,y)\in X\times Y$ define 
\[ \tilde  g^{b}_{x,y} = \left\{
 \begin{array}{lll}
 1       & \spacenine    & g^{b}_{x,y} \geq \overline\ep \\
 g^{b}_{x,y} / \overline\ep & & g^{b}_{x,y} < \overline\ep.
 \end{array}
 \right. \]
Define the {\em auxiliary absorption rate} by
$ \tilde  a_{x,y}(s) = 
 \sum_{b \in B} y^s_b \tilde  g^{b}_{x,y}. $
Note that
$a(s) \leq \tilde  a(s) \leq {a(s)/ \overline\ep}$. 
$$\mbox { Define }
\tilde  v^2(b) = (1 - {g^{b} \over \tilde   g^{b}}) 
r^2(s) + {g^{b} \over \tilde   g^{b}} v^2(b)\spacenine \spacenine (6)$$
with $\tilde   v^2(b) :=r^2(s)$ if $g^{b}=\tilde   g^{b}=0$.

Next we need to use large  quantities $Q_1>1$ and $Q_2>1$ that will 
 be determined precisely later (in the next section) 
 by the choice of  $\alpha$, $\ep$,
 $\overline \epsilon$, $\hat \ep$ and    
$\tilde \epsilon$.
        Define $L:= Q_1Q_2$ and 
define   $K:=L^{|{\cal N}|}$.

Define the function $h : [1,\infty) \to [1,K]$ by
$ h(r) =  \min\{r, K\}$. Order the members $\{ s_1, \dots , s_m\}$ 
 of ${\cal N}$ with
 $\tilde  a_{x,y}(s_1) \leq \tilde  a_{x,y}(s_2)\leq 
 \dots \tilde  a_{x,y}(s_m)$. 
Define 
\[ \tilde  w_{x,y}(s_k) = 
\prod_{j=k}^{m-1}
  h( {\tilde  a_{x,y}(s_{j+1}) 
\over \tilde  a_{x,y}(s_{j})}). \]

  For any move 
 $b$ at a state $s\in {\cal N}$ define  
$\overline g_{x,y}^{b}$ to satisfy  
$$(1-\tilde  g^{b}_{x,y})
 = (1-\overline g^{b}_{x,y}) (1-g^{b}_{x,y}). 
  \spacenine \hfill (7)$$
If $\tilde  g^{b}_{x,y} = 1$,
 then $\overline g^{b}_{x,y} = 1$ as well.
Note that 
$$g^{b} v^2(b) + (1-g^{b}) \overline g^{b} r^2(s) = 
g^{b} v^2({b}) + (\tilde  g^{b} - g^{b}) r^2(s) = 
\tilde  g^{b} \tilde  v^2(b). 
 \spacenine \hfill (8)$$

 For every $s\in 
{\cal N}$ and $h\in \tilde {\cal H}$ denote
$ N^s(h) = \#\{ n \in {\bf N} \ |\ s_n = s \} \in {\bf N} \cup \infty. $
 For $1 \leq i \leq N^s(h)$
 let $n^s_i(h)$ be the stage with 
 the $i$th occurrence of the state $s$ in $h$.
 If the initial state of $h$ is $s$, then $n^s_1 = 0$ and $N^s (h)\geq 1$.

 Define the discounted 
 evaluation of a move $b$ at a state $s\in {\cal N}$ according to 
$$\xi^{b}_{x,y} = E^b_{x,y} r^2(h) \ [\quad 
\sum_{i=1}^{N^s(h) - 1} 
\overline g^{b_{n^s_i(h)}}(1-{\delta\over \tilde  w_{x,y}(s)})^{i-1}
\prod_{k=1}^{i-1} (1-\overline g^{b_{n^s_k(h)}})+ $$  
$$(1-{\delta\over \tilde  w_{x,y}(s)})^{N^s(h) - 1} 
\prod_{k=1}^{N^s(h) - 1} (1-\overline g^{b_{n^s_k(h)}})\quad ]
\ , \spacenine (9)$$
 where $E^b_{x,y}$ stands for the expectation
 over all infinite 
histories $h\in \tilde {\cal H}$ with initial state $s_0 = s$,
 assuming that Player Two plays the action $b$ at stage $0$, the first stage,
and afterwards follows $y$, whereas Player One follows $x$ always.
 \vskip.2cm

{\bf Lemma 5.1} The function $\xi^b_{x,y}$ obeys the properties
$$\xi^{b}_{x,y} = 
\tilde  g^{b}_{x,y} \tilde  v^2_{x,y}(b) + 
(1-{\delta\over \tilde  w_{x,y}(s)}) (1-\tilde  g^{b}_{x,y}) 
\xi_{x,y}(s) \spacenine \hfill (10)$$
and
$$r_{x,y}^2(s) = 
\xi_{x,y}(s) 
\left( 1 + {\delta ( 1 - \tilde  a_{x,y}(s) ) \over 
        \tilde  w_{x,y}(s) \tilde  a_{x,y}(s) } \right)
\spacenine \spacenine \hfill (11) $$
where $\xi_{x,y}(s) = 
\sum_{b \in A^s_2} y_b^s \xi^{b}_{x,y}$.

\vskip.2cm {\bf Proof:}

We now verify that $\xi$ satisfies  (10) and (11).
Separate the summation in (9) into three parts.
\begin{itemize}
\item   All histories such that $N^s(h) = 1$.
        The probability of this event is $g^{b}$, and 
        the conditional expectation is $v^2(b)$.
\item   All histories such that $N^s(h) > 1$ and $i = 1$.
        The probability of this event is $1-g^{b}$, and
        the conditional expectation is $\overline g^{b} r^2(s)$.
\item   All histories such that $N^s(h) > 1$ and $i > 1$.
        The probability of this event is $1-g^{b}$.
        Factor out one power of
 $(1-\overline g^{b})(1-{\delta\over \tilde  w})$;
        the conditional expectation is
        $(1-\overline g^{b})(1-{\delta\over \tilde  w})\xi(s)$.
        By (7) this part contributes 
        $(1-\tilde  g^{b})(1-{\delta\over \tilde  w})\xi(s)$ to the sum.
\end{itemize}
Putting together the three parts, with (8) connecting the first two parts,
 we get (10).  For equation (11) we use (10) and 
 take the expectation with respect to the moves. 
\hfill $\Box$ \vskip.2cm 

Notice that formula (11) is a slight variation of the standard 
 relationship between discounted and undiscounted evaluations. 
$\xi$ will serve as the auxiliary discounted payoff evaluation of Player Two.
Note that
$r_{x,y}^2(s) \geq \xi_{x,y}(s)
 \ \forall s \in {\cal N}.$ 
 Define $\overline \xi_{x,y}(s)$ to be maximal value 
 $\max_{b\in A^s_2} \xi_{x,y}^b$.    
\vskip.2cm

{\bf Lemma 5.2:} 
For every $s,t \in {\cal N}$, $\gamma>0$, and 
$(\delta,x,y) \in (0,1] \times X \times Y$   
\begin{itemize}

\item   $  \tilde a_{x,y}(t) \leq
                K \tilde  a_{x,y}(s)$ 
        \mbox {implies that}  
        $\tilde  w_{x,y}(t) \tilde  a_{x,y}(t)
                \leq \tilde  w_{x,y}(s) \tilde  a_{x,y}(s)$,  
\item    $  \tilde  w_{x,y}(s) \tilde  a_{x,y}(s)
                \leq 
 \tilde  w_{x,y}(t) \tilde  a_{x,y}(t)$ and 
 $r_{x,y}^2(s) \leq r_{x,y}^2(t)+\gamma$ imply 
 that  $\xi_{x,y} (s)\leq \xi_{x,y} (t)+ \gamma +\delta$.
 \end{itemize}
 \vskip.2cm 
{\bf Proof:} 
The first  part follows directly from the definition of $\tilde w$. 
For the second part,
note that for every $r,\tilde w,a > 0$ and $0 < \delta < 1$
\[ {r\tilde wa \over \tilde wa
 + \delta(1-a)} - {r\tilde wa \over \tilde wa + \delta} =
{r\tilde wa^2 \delta \over 
(\tilde wa + \delta)(\tilde wa + \delta - \delta a)} \leq
{r \tilde w a^2\delta^2 \over
 \tilde w^2 a^2} = {r\delta \over \tilde w}. \]
Moreover, $r\tilde wa\over \tilde wa +
 \delta$ is an increasing function in $\tilde wa$.
  Given $\tilde w\geq 1$, from the above  we have that 
                  $r^2(s)$ and $\tildew(s) \tilde a(s)$  determine 
       $\xi(s)$ except for a quantity of no more than $\delta$. 
     \hfill $\Box$ \vskip.2cm


\subsection{The Best Reply Correspondence}

For every state $s \in {\cal N}$  define 
\begin{eqnarray}
B_{\delta,1}^s(x,y)
 &=& \argmax_{a \in A^s_1}\ w^1_{x,y}(a) 
        \spacenine  \nonumber \\
B_{\delta ,2}^s(x,y)
 &=& \argmax_{b \in A^s_2} \ \xi^{b}_{x,y}
\hbox{ if }\xi_{x,y}(s) > j^{\alpha}_{x}(s) \nonumber \\ 
B_{\delta ,2}^s(x,y)
 &=& J_x^{\alpha}(s)\cup \argmax_{b \in A^s_2} \ \xi^{b}_{x,y}  
\hbox{ if }\xi_{x,y}(s) = j^{\alpha}_{x}(s). \nonumber \\
B_{\delta ,2}^s(x,y)
 &=& J_x^{\alpha}(s)
\hbox{ if }\xi_{x,y}(s) < j^{\alpha}_{x}(s). \nonumber
\end{eqnarray}
 Player One maximizes her un-discounted payoff, while Player Two
maximizes his auxiliary payoff, given 
 that it is not too small. 

Let the corresondences $\overline B_{\delta,1}^s$
  and $\overline B_{\delta ,2}^s$ be those  defined 
 by the closure of the graphs of the correspondences 
 $B_{\delta ,1}^s$ and  $B_{\delta ,2}^s$ 
in $(X\times Y)\times A^s_1$ and 
 $(X\times Y)\times A^s_2$, respectively. 
Define conv $(\overline B_{\delta,1}^s)$ and 
 conv $(\overline B_{\delta,2}^s)$ to be 
 the correspondences with  graphs in $(X\times Y)\times X^s$ and 
 $(X\times Y)\times Y^s$, respectively, such that 
 $z\in $ conv $(\overline B_{\delta,1}^s(x,y))$ if and only if  
 $\{ a\in A^s_1\ | \ z_a >0\}$ is a subset of
 $\overline B_{\delta ,1}^s(x,y)$ and 
 $z\in $ conv $(\overline B_{\delta,2}^s(x,y))$ if and only if  
 $\{ b\in A^s_2\ | \ z_b >0\}$ is a subset of
 $\overline B_{\delta ,2}^s(x,y)$.
Define the correspondences $B_{\delta ,1}$  from 
 $X\times Y$ to $X$ so that $(x,y)$ in the domain corresponds 
 to the sets $B_{\delta,1}^s (x,y)$ in the range, and likewise 
 define the correspondence $B_{\delta, 2}$ from 
 $X\times Y$ to $Y$.
We define the correspondence $F_{\delta} : X\times Y \rightarrow 
 \rightarrow X\times Y$ by 
 $F_{\delta} (x,y) = (B_{\delta ,1}(x,y), B_{\delta ,2}(x,y))$. 
By Kakutani's fixed point theorem for every $\delta>0$ the correspondence 
$F_{\delta}$ has a fixed point.

\subsection{Two Lemmas on Fixed Points}

     We assume in the rest of the section that $(x,y)$ is 
      a fixed point for $F_{\delta}$. We prove Lemmatta 5.4 and 
 5.5, described in the introduction. 
\vskip.2cm 

{\bf Remark 5.3:} 
Since the jump correspondence
 is used before $\xi$ gets close to $0$, any fixed 
 point $(x,y)$ of $F_{\delta}$ is absorbing.  This implies that 
$r_{x,y}^2(s) \geq j^{\alpha}_{x}(s) \quad \forall s \in {\cal N}.$ 
Indeed, suppose for the sake of contradiction 
that $r_{x,y}^2(s) < j^{\alpha}_{x}(s)$. 
Denote by $e$ the stopping time that is defined by the
first stage in which the game leaves the set 
$\{u \ |\ \xi_{x,y}(u) < j^{\alpha}_{x}(u)\}$.
Recall from Section 2 that $j^{\alpha}_x$ is a sub-martingale. 
Since   for every absorbing state $s\in {\cal A}$
$\xi(s) = j^{\alpha}(s)=r^2(s)$ 
we have $ j^{\alpha}_{x}(s) \leq \E j^{\alpha}_{x}(s_e) \leq 
        \E \xi_{x,y}(s_e) 
        \leq \E r^2_{x,y}(s_e) = r^2_{x,y}(s), $
as desired.
\vskip.2cm

{\bf Lemma 5.4} If $\overline \ep \leq \omega \alpha /4$ then 
there is a choice for  $L^*>1$ and $\delta^*>0$ such that if $L\geq L^*$ 
 and $0<\delta \leq  \delta^ *$  and 
  $(x,y)$ is a fixed point of $F_{\delta}$  then \newline 
 1) $\xi_{x,y}(s) \geq j^{\alpha}_{x}(s)$ for all $s\in {\cal N}$, \newline 2) 
 if the jump correspondence is used at $s$ 
then  $\xi_{x,y} (s) 
\leq r_{x,y}^2(s) -3\overline\ep$\newline 3) for any  action $b$ from
 $J_x^{\alpha}(s)$ 
 used in $y^s$   
$g^{b}_{x,y} < \overline\ep$ , 
and \newline 
4) the overall probability that Player Two plays an action from 
$J_x^{\alpha}(s)$
at any $s\in {\cal N}$ is at most $\omega \alpha/20$.
\vskip.2cm

{\bf Proof:}
Let $L^*= {100 |{\cal N}| \over \omega^2 \alpha ^2 \overline \ep}$ and 
 $\delta^* = \overline \ep \alpha ^3 \omega ^3/(300 |{\cal N}|)$.  
Choose $t$ to be a member of 
 ${\cal N}$ with the largest difference 
 $j^{\alpha}_x(t)- \xi (t)$, and we must presume that 
 this difference is non-negative.  
We will show that this difference can be no more than $0$ and that 
the    frequency devoted to the jump correspondence at any such 
 state can be no 
 more than $\alpha \omega/20$.

We presume for the sake of contradiction that the frequency devoted 
 to the jump correspondence at $t$ is at least $\alpha \omega /20$. 
  Since $r^2\geq j^{\alpha}_x$ the 
 expected value of the jump function $j^{\alpha}_x$ at  
  the states reached on the next stage after $t$ using 
the jump correspondence $J^{\alpha}_x$ is 
 at least $\alpha \omega$ more than  $j^{\alpha}_x(t)$, 
 we must assume for any move  from  $J^{\alpha}_x(t)$ that 
there is at least one state $u$ reached by this move  
  with a probability of at least $\alpha^2 \omega^2
 \over 40 |{\cal N}|$ such  that
 $j_x^{\alpha}(u) \geq j_x^{\alpha}(t)+ \alpha \omega /2$, 
 necessarily with $\esc (u,t)\leq \alpha \omega /4$. (If $\esc (u,t)> \alpha 
 \omega /4$ then   
 $ a (t)\geq \alpha^3 \omega^3 /(160 |{\cal N}|)$ and 
 by (11) and the size of $\delta ^*$ we have made $\xi(t)$ too 
 close to $r^2(t)$ contradicting  $j_x(t)\leq r^2(t)-\alpha \omega/2$, --
 which must follow by Remark 5.3 since otherwise any move from 
 the jump correspondence would be evaluated in an undiscounted 
 way strictly above the level $j^{\alpha}_x(t)$.) 
  By the definition 
 of $\tilde w$, the size of $L^*$  and (3)
 we have $\tilde w (t) \tilde a(t) \geq \tilde w (u) \tilde a(u)$. 
By $\esc (u,t)\leq \omega \alpha /4$ it follows 
 that $|r^2(t)-r^2(u)|\leq \omega \alpha /4$.
  But by Lemma 5.2 we have $\xi (t) \geq 
\xi (u)- \delta -\alpha \omega /4 $.
  With  the size of $\delta^*$  this 
 contradicts $j_x^{\alpha}(u)\geq j_x^{\alpha}(t) + \alpha \omega /2$ and 
 the choice of $t$.

             Next, suppose for the sake of contradiction that  
         $J^{\alpha}_x$ is used at $s$ and 
    $g^{b}\geq \overline\ep$ for some move   
             $b \in J^{\alpha}_x(s)$. 
Indeed,  $g^{b} \geq \overline\ep$ implies that $\tilde  g^{b} = 1$.
In particular,              using Remark 5.3, 
$ \xi_{x,y}^{b} =  
w_{x,y} (b) \geq \sum_t p(t|s; x,b)  j^{\alpha}_{x}(t) \geq 
j^{\alpha}_{x}(s) + \omega\alpha.$
Thus, for every $b' \in B^2_\delta(x,y)$ that maximizes $\xi$,
$\xi ^{b'}_{x,y}\geq 
\xi ^{b}_{x,y}\geq r^2_{x,y}(s)  \geq 
j_{x}(s) + \omega\alpha/2$.
Since the overall probability to play actions from the jump 
 correspondence   is smaller
than $\omega \alpha/20$, 
this contradicts the assumption $\xi(s)\leq j^{\alpha}_x (s)$.

     Now we presume for the sake of contradiction  
     that $\xi(s) \geq r^2(s)-3\overline\ep$ and the $J^{\alpha}_x$
 correspondence is       used at $s$.
      Since we must assume that $\xi(s)=j^{\alpha}_x(s)$, we have 
      an increase in the value of $r^2$ of at least 
       $\omega \alpha-3\overline\ep$ from a move in $J_x^s$.  
       By the dominance of $\omega \alpha$ over 
        $4\overline\ep$, we must conclude that $g^b>\overline \ep$, a 
          contradiction.
              \hfill $\Box$
                                           \vskip.2cm

Lemma 5.4 is the most problematic aspect of extending this proof 
 to the case of finitely many positions. Any identification of infinitely 
 many states as a single state may be meaningless if the states 
 reached from it are not also identified. A more flexible definition
 of the discounted evaluation may be necessary. 
For example, at a state $s$ one could discount future visits to 
 other states $t$ according to the difference between  Player Two's 
 undiscounted  expected payoffs from these two states.

The following lemma claims that if the auxiliary payoff is too far from the
real payoff and  the action causes absorption with small probability,
then this probability is very small. This radical discontinuity 
 is the key argument to our whole approach. 
 \vskip.2cm 
{\bf Lemma 5.5}                For  $L$, $\alpha$, $\overline \ep$
 and $\delta$ satisfying the conditions of 
   Lemma  5.4 and $(x,y)$ a fixed point of $F_{\delta}$  
 if $\xi (s) \leq r^2(s) - 2\overline\ep$ and 
$g^{b} \leq \overline\ep$
then
$g^{b} \leq 1.1\ \delta \xi (s) / \tilde   w (s)$
 and $g^{b}\leq 2.3 \ \overline\ep \tilde   a(s)\leq 2.3 \  a(s)$. 
\vskip.2cm

 {\bf Proof:}
 First we claim that 
$\overline\xi (s) - \xi (s) 
< { \delta\alpha \omega  \over 19 \tilde   w (s)}\xi(s)$. 

                                         If the jump
 correspondence  at $s$ is used 
 and   $b$ is such a move,
since $g^{b} \leq \overline\ep$
 (from Lemma 5.4) it follows that $\tilde  g^{b} = g^{b} 
/ \overline\ep$.
Hence from (6) we have 
$$\tilde v^2(b)  = 
(1-\overline\ep) r^2 (s) + 
\overline \ep v^2(b)  
 \geq r^2 (s) - \overline\ep 
\geq \xi  (s) + \overline\ep.
\spacenine \hfill (12)$$
Moreover, from (10) and (12) we have
$$\xi^{b}  \geq  \xi (s)+
\tilde   g^{b}  (\tilde   v^{b} -\xi (s)) 
- \delta \xi (s)/\tilde   w (s) \geq 
  \xi (s)(1
- \delta/\tilde   w (s)), $$ 
 and by Lemma 5.4, since $\xi (s)$ is 
 the average of $\overline \xi (s)$ and such $\xi ^{b}$, 
 we have $(1-\alpha \omega /20) (\overline \xi (s)- 
 \xi (s))\leq {\delta \alpha \omega \xi (s) \over 
 20 \tilde   w (s) }$, so the claim follows.

                         Considering now any move $b\in A^s_2$ 
                          that is used with 
              $g^{b}\leq \overline\ep$ and  
 looking again at formula (10) we have 
 $\overline \xi (s)\geq \xi ^{b} \geq 
 \tilde   g^{b} (\tilde   v^{b} -\xi (s)) 
  + (1-\delta/ \tilde   w(s))\xi(s)$ and hence  
 $\tilde   g^{b}(\tilde   v^{b}-\xi(s)) 
  \leq 1.1\ \delta 
  \xi (s)/\tilde   w(s)$, since by the above claim  
   $\overline \xi (s)- \xi(s)$ is  small 
    compared to    
    ${\delta\over \tilde   w(s)}\xi(s)$. 
   First consider the consequence of 
   $\tilde   v^{b}-\xi (s)\geq 
    \overline\ep $, namely $ g^{b}=\overline\ep \tilde   g^{b}
  \leq 1.1 \ \delta \xi (s)/\tilde   w(s)$.
Second, consider  
$\tilde   g^{b}\leq {1.1\ \delta  \xi(s)
\over (\tilde   v^{b}-\xi(s))\tilde   w(s)}\leq 
 {1.1\ \delta  \xi(s)
\over (r^2(s)-\xi(s)-\overline\ep)\tilde   w(s)}$, proven above.  
 Since ${r^2(s)-\xi(s)-\overline\ep \over  r^2(s)-\xi(s)} \geq 1/2$, 
  we get 
$\tilde   g^{b}\leq 
 {2.2 \delta  \xi(s)
\over (r^2(s)-\xi(s))\tilde   w(s)}.$   
 Now apply formula (11) for 
$2\overline\ep \leq r^2(s)-\xi(s) = \xi(s){\delta (1-\tilde   a(s)) \over  
 \tilde   a(s)\tilde   w(s)}$.   
          Since     $\tilde   w(s)\geq 1$ and    $\xi(s)\leq 1$ 
           we have $\delta (1-\tilde   a(s)) 
\geq 2\overline\ep \tilde   a(s)$, and 
            from $\delta < 
 \overline\ep /25$ we have $\tilde   a(s) \leq 1/50$. This allows  
  us to conclude with 
$ {g^{b}\over \overline \ep} = \tilde   g^{b}\leq 
 {25\over 24} 2.2\ \tilde   a(s)\leq 2.3 \  
 \tilde   a(s)\leq {2.3 \   a(s)\over \overline \ep}$.   
\hfill $\Box$ \vskip.2cm


\section{Second Main Theorem}

The goal of this section is to prove Theorem 2, which states 
  that the conditions of Theorem 1 are always satisfied. 
 First we need a simple but useful lemma. 

 \vskip.2cm 
{\bf Lemma 6.1}
 For every two distinct non-absorbing 
 states $s,t$  with $ \esc (t,s)\leq 
 \gamma<1$ in an absorbing time homogeneous Markov chain
$P^t(t,s) \mu(s,t)$ does not differ from $ a(t)$ by more than 
 a factor of $2\gamma$,  
 $\esc(t,s) / \mu(t,s)$ is within a factor of $3\gamma$ to the 
 ratio  that, starting at  $t$ or $s$,
 the last visit before absorption was at $t$ rather than 
  at $s$.  Furthermore, with or  without 
 the assumption that the Markov chain is absorbing and with a 
 start at either $s$ or $t$, the ratio of 
 the expected number of visits to $s$ to those at $t$ is at least 
 $1-4\gamma$ times  the ratio of $P^t(t,s)$ to $P^s(s,t)$.

\vskip.2cm {\bf Proof:}
The first two claims   follow directly from the formulas (1) and (2). 
 The third claim follows from the first claim if the Markov chain 
 is absorbing.  Otherwise we recognize in $1/P^t(t,s)$ the expected 
 number of visits to $t$ before reaching $s$. 
\hfill $\Box$ \vskip.2cm 
   
{\bf Remark 6.2}
At a fixed point of $F_{\delta}$ satisfying the properties of 
 Lemmatta 5.4 and 5.5,
 if $\xi(s)\leq r^2(s)-2\overline\ep$ and $b$ is a move 
 at $s$ with  
 $g^b\geq \overline \ep$ then $w^2(b)= \overline \xi(s)$, 
 which is by Lemma 5.5 also within $\delta/20$ of $\xi(s)$.   
 \vskip.2cm 

{\bf Theorem 2:} For any choice of positive 
 $\ep$, $\overline \ep$, $\hat \ep$, and 
 $\tilde \ep$ satisfying the inequalities stated in Theorem 1  
  all  conditions  of Theorem 1 are satisfied.   \vskip.2cm  

{\bf Proof:} 
Because it is sufficient to demonstrate the conclusion of Theorem 1 
 with smaller choices for $\overline \ep$, $\hat \ep$ and  
 $\tilde \ep$, we will assume without loss of generality that $\alpha$ 
 is small enough so that for every $s\in {\cal S}$ $c^{\alpha}(s)$ is 
 within $\ep /2$ of the undiscounted zero-sum value $c_2(s)$, as described 
 in Section 2, and $\overline \ep < \alpha \omega /4$.

  Define $\beta := 
{1\over 2} \tilde \ep \hat \ep^{|{\cal N}|}/(3^{|{\cal N}|}|{\cal N}|)$. 
 We require that   $L:=Q_1Q_2$   is large enough to
 satisfy the   conditions of Lemma 5.4  and also that\newline    
 $Q_1>   80 |{\cal N}|^3  m^2 / (\rho \ \overline \ep ^2 \
\hat \epsilon^2 \ \tilde \ep ^2 \ \beta^2)$
 and $Q_2> 80m  |{\cal N}| / 
(\rho \ \overline \ep \ \hat \epsilon \ \tilde \ep)$.

We begin with $\delta$ sufficiently small, so that the 
 condition of Lemma 5.4 holds. 
Next, we consider fixed points of $F_{\delta}$
 corresponding to decreasing $\delta>0$ 
 that have convergent subsequences for certain variables living in 
 compact spaces --  
 the stationary strategies 
 in the space $X\times Y$, the values 
 $\nu (a,b)$  for all pairs of moves at all states, the expected payoffs 
 $r^1(s)$, $r^2(s)$,    and the absorption 
 rate $a(s)$ for every $s\in {\cal N}$, and 
 the probabilities $\esc (t,s)$ for all pairs of states.

We define a move $a\in A^s_1$ or $b\in A^s_2$ to 
 be a {\em limit} move if and only if the frequency of its use 
 does not converge to zero as $\delta$ goes to zero, and define  
 $\hat q$ to be the minimal positive limit value 
 for a frequency of a limit move  chosen by either player.
We define the quantity $\hat \mu$ to be the minimal positive limit value 
 for $\esc (s,t)$,  $\hat \nu$ 
 to be the minimal positive limit value for $\nu (a,b)$, and 
 $\hat a$ the minimal positive limit value for $a(s)$.

Next we must define  
   the partition ${\cal R}$ of a subset $P$.  Define 
 a directed graph  on the  space ${\cal N}$  
 by $t\rightarrow s$ if and only if in the limit 
 $\esc (t,s)$ approaches zero.
  The relation is transitive, but not necessarily 
 symmetric.
      It has an additional property, that if $t\rightarrow s_1$ 
 and $t\rightarrow s_2$ then either $s_1 \rightarrow s_2$ 
 or $s_2 \rightarrow s_1$.   This is easy to confirm, because if 
 $s_1$ was not reached with probability approaching 
 one  on the way from $t$ to $s_2$ 
 then it must be reached with probability approaching one
 after the state $s_2$.  
 Next define a relation $\sim$ that 
 is  symmetric, transitive, and reflexive on a appropriate subset; 
 $s\sim t$ if and only if $\mu (s,t)$ approaches zero, and 
 $s\sim s$ if and only if $a(s)$ approaches zero.  $\sim$ defines 
 a partition ${\cal P}$ of a subset $P'$ of ${\cal N}$.  Now 
 we relate $\rightarrow$ to $\sim$.  
Define ${\cal R}$ to be the subset 
 of ${\cal P}$ defined by $A\in {\cal R}\subseteq {\cal P}$ if and only 
 if  $u\in A$ and $u\rightarrow s$ implies that $s\in A$.   
 Any state $s\not\in A\in {\cal R}$
 such that $\esc (s,u)$ approaches zero for any (equivalently 
 some) state $u\in A\in {\cal R}$ is called a {\em satelite} of $A$. 
 Due to the above, a satelite  of
 $A\in {\cal R}$ cannot be a satelite of any other member 
 of ${\cal R}$ and every member of $ Q\in {\cal P}$ such that 
 $Q$ is not in ${\cal R}$ must be a satelite of the same  
 $A\in {\cal R}$. We call an primitive exit $(a,b)$ from 
 $R\in {\cal R}$ to be a {\em satelite exit} if with certainty 
 the exit results in motion that doesn't leave   $R$
 or its collection of  satelites. 

For every $R\in {\cal R}$ we define the set $B_R$ of Player Two moves 
 in $R$ to be $B_R:= \{ b\in A^s_2\ | \ $for some limit move 
 $a\in A^s_1$ $(a,b)$ is an primitive exit 
 from $R$ that is not a satelite exit$\}$. 

If $s$ is a satelite of $R\in {\cal R}$ then 
in the limit the probability  that the last 
 visit to the pair $s$ or  any $u\in R$  was the state $s$  
must go to zero. Therefore     
 $\nu (a,b)$ approaches zero for 
 any satelite exit  $(a,b)$ at $u\in R\in {\cal R}$.
These facts follow directly from Lemma 6.1 and 
 $\esc (s,u) / \mu (u,s)$ going to zero 
 in the limit.

We show  for every $R\in {\cal R}$ and pair 
 $s,t\in R$ that the probability of using some exit  
 in $E^{B_R}(R)$ before reaching $t$ from $s$ also 
 approaches zero.  First this holds for 
 any non-satelite primitive exit
 from $R$, because the probability 
 of reaching a non-satelite outside of $R$  
 would be at least $\rho$ and therefore the probability 
 of absorbing before  reaching $t$  must be in the limit at least 
 the probability of using this exit  times $\hat \mu \rho$. 
   The same 
 arguments holds for the use of any move in $B_R$, but with the 
 quantity  
  $\hat \mu \rho \hat q$ instead of $\hat \mu \rho$. 

 More difficult is to  show that the above holds 
 for any  satelite exit $(a,b)$ at $u\in R$.  Let $v$ be any 
 satelite of $R$ reached with positive probability from this 
 exit.  Let $\pi$ be the probability of using 
 $(a,b)$ before reaching $t$ from a start at $s$ and let $\theta$ be 
 a bound on the probability of not reaching any member 
 of $R$ from any other member of $R$ or from a satelite 
 of $R$.  Let $\hat \gamma$ be the probability of reaching $v$ from $s$ before 
 reaching $t$, with  $\hat \gamma \geq \pi \rho$. Going through 
 the state $v$, the probability of reaching $t$ is 
 at least $1-\theta$ and the combined probability of reaching 
 $t$ from $s$ is also at least $1-\theta$. This means that 
 the probability of  reaching $t$ from $s$ conditioned 
 on not going through $v$ is at least  $1-{\theta\over (1-\hat \gamma)}$. 
 So conditioned on not arriving at $v$ before $t$ 
there is 
 at most a $2\theta\over 1-\hat \gamma $
 probability of absorbing before getting 
 back to $s$. In the limit $2\theta \over 1-\hat \gamma$ cannot stay 
 above $1$, because $\theta$ goes to zero and in the limit 
 $\hat \gamma$ cannot go above $1-\hat \mu$. 
This means that eventually
 the probability of reaching $v$ from $s$ must be at least 
 $\hat \gamma \sum_{i=0}^{\infty} 
(1-\hat \gamma)^i(1- {2\theta\over 1-\hat \gamma})^i= 
 \hat \gamma \sum_{i=0}^{\infty} (1-2\theta-\hat \gamma)^i=
 {\hat \gamma  \over 2\theta+ \hat \gamma}$. 
 But this  probability 
 to reach $v$ from $s$ cannot go  
 above $1-\hat \mu$ in the limit, which is possible only if 
 $\pi$ goes to zero as $\theta$ goes to zero also.

   Define $\epsilon^*$ to be $( \hat \nu \hat \mu \hat q \hat a 
 /  K) ^{3|{\cal N}|}$. 
We require of a fixed point of $F_{\delta}$ 
that the values for which we have convergent subsequences 
 are within  $\epsilon ^*$ of their limit values. 
 We require that the probability of using 
 any  exit before moving  from any $s$ to $t$ for any 
 pair  $s,t\in R\in {\cal R}$ 
  is no more than $\ep^*$ and 
  for every $R\in {\cal R}$ that the sum of  $\nu (a,b)$ over 
 all  the satelite 
 exits $(a,b)\in E^{B_R}(R)$ is no more than $\ep ^*$ (as demonstrated 
 above). Furthermore 
 we require that $\delta < (\ep^*)^2$. We let 
 $(x,y)$ be a fixed point of $F_{\delta}$ satisfying these properties. 
 If the stationary strategy is not specified, then $(x,y)$ is 
 intended.  
\vskip.2cm

  {\bf Step 1; For every $s\in R\in {\cal R}$ show that if 
 $z_{x,y}^{2.5\overline \ep}(s)\geq \beta $ 
then there exists an $\tilde \ep$ simplication $\overline x$ of $x$ 
  such that 
 $z^{2.4\overline  \ep}
_{\overline x, y}(R)\geq 1- 3|R|/Q_1\beta$ and 
 for all $t\in R$ that  
$z^{2.4 \overline \ep}
_{ x, y}(t)/ z^{2.4 \overline  \ep}_{x,y}(s) \geq 
 (1-{m\over 2 \rho Q_2 \overline \ep} -{4m|R|\over \rho Q_2}) 
(z^{2.4\overline  \ep}
_{\overline x, y}(t)/ z^{2.4\overline  \ep}
_{ \overline x, y}(s))$:} \vskip.2cm

 For every $d \geq 1$ define
\[ T_d = \{ t \in R \ |\ \mu(s,t) \leq d\tilde  a(s)\} \cup \{s\}. \]
Denote $T = T_d$, where 
$d \in (1 ,L^{|{\cal N}|-1})$ satisfies
$T_{Ld} \setminus T_d = \emptyset$.
Since $K = L^{|{\cal N}|}$, for every $t \in T$ we have 
$\tilde  a(t) \leq a(t)/\overline \ep \leq \mu(s,t)/
\overline \ep \leq d\tilde  a(s)/\overline \ep \leq 
        K\tilde  a(s)$, and it follows that 
$ \tilde  w(t) \tilde  a(t) \leq \tilde  w(s) \tilde  a(s). $ 
 With $\xi(s)\leq r(s)-2.5\overline \ep$ we have by (11) 
 that $\tilde a(s)\leq \delta /\overline \ep$, 
 $a(s)\leq \delta /\overline \ep$,
 and $\mu(s,t)\leq \delta K/ \overline \ep < \ep ^*$, meaning 
 that $T$ is a subset of  $R$. 
Since $t\in T$ satisfies  $|r^2(s)- r^2(t)|\leq \epsilon^*$ we have  
$\xi(t) \leq \xi(s) +(\epsilon^* +\delta) $.

Define a  quantity 
\[ p^t = \left\{
\begin{array}{lll}
a(t) / Q_2\mu(s,t)        & \spacenine    & t \in T \setminus \{s\} \\
1/Q_2                     &               & t = s
\end{array} \right. \]
Define the stationary strategy $\overline x $ by removing 
  from $x$ all  Player One moves at states $t\in T$ that are played with
probability 
smaller than $p^t / \rho$, and normalize the remaining vector. This means 
 that if $u$ is reached in one stage from $t\in T$ by $\overline x$ and 
 a Player Two move $b$, then $p(u|t; \overline x,b) \geq p^t$.

          We use critically       from Lemma 6.1
          that $a(u)/\mu (u,s)$ is approximately $ P^u(u,s)$ (within 
 a factor of $2\epsilon ^*$)
 for any $u\in R$, so that 
from    Lemma 3.2 and Lemma 6.1 with the change from $(x,y)$ 
 to $(\overline x,y)$  the ratio of visits  
       at  $t\in R$ to those at $s$ 
       cannot increase by  more than a  factor of ${8|T|m\over \rho Q_2}$. 
 Furthermore, by the definition of $\overline x$, $\hat \nu$, 
 $\delta \leq (\ep^*)^2$ and Lemma 5.5 
  there are no non-satelite exits performed inside of $T$ 
  other than those generated by Player 
 Two moves $b\in A_2^t$ with   
  $w_{x,y}(b)=\xi_{x,y} (t)$.  
        Combined with the fact  
        that the absorption
 rate of any move $b$ with $g^b\geq 2.4 \overline \ep$
 is  altered by a factor or no more than $m/(2\rho \overline \ep Q_2)$ by   
         the switch to $\overline x$ and  that  $2\ep^*$ is greater 
 than the   
       probability that the last visit to $R$ was at a  
  satelite exit,  
we  have everything but the claim that there is only insignificant 
 motion  with $(\overline x,y)$ toward absorption from states in $R$ 
 outside of the set $T$.

    We can break up the 
       absorption from $R$ generated by the strategies $(x,y)$ 
       in terms of where was the last visit  in  $R$.    
Let $t \in T$, $u \in R \setminus T$ and $b \in B$ be a move 
such that
$p(u|t,\overline x,b) > 0$, necessarily with $p(u|t; \overline x,b)\geq p^t$.
To complete the claim of Step 1  
  it suffices to show that
        ${\esc_{x,y} (u,t)\over \mu_{x,y} (u,t) }
        \leq {2.5 
\over   Q_1}$ for every such $u\in R\backslash T$.
\vskip.2cm 

  {\bf Case 1; $u\in R \backslash T$ is reachable from  $ T$ \underline 
  {only} 
   by  Player Two  moves $b$ with 
 $g^b\geq \overline \ep$:}
      
It follows immediately from the fact that Player Two has no more than 
 $m|R|$ moves in $R$   
that  ${esc} (u,t)\over \mu (u,t)$ is smaller than $2.5/Q_1$, since 
 any such move doesn't return to $R$ with a probability of at least 
 $2.5\overline \ep$ and  with  
 at least $1-2\ep ^*$ probability there is motion 
   from $u$ back to   
  $t\in T$. 
                               \vskip.2cm 

{\bf Case 2;  $t\not=s$, and 
   $u\in R \backslash T$ is reachable by $(\overline x, y)$  
   from a $t\in T$  
   by a  move $b$ of Player Two with $g^b< \overline \ep$:}

By Lemma 5.5 we have 
$$p^t \esc(u,t) \leq g^{b} \leq 
 2.3 \  a(t).$$ 
 Since $p^t = {a(t) \over  Q_2 \mu(s,t)}$ we have
$\esc(u,t) \leq 2.3 \  Q_2 \mu (s,t)$. Since $\mu$ is a metric we have   
 from 
 $\mu (s,u)\geq L\mu (s,t)$
$$ {\esc(u,t) \over \mu(u,t)} \leq {2.3 \ Q_2 \mu (s,t) \over \mu (u,t)}\leq
{2.3 \ Q_2 \over L-1} \leq {2.4\over Q_1}.$$

{\bf Case 3; 
    $u\in R \backslash T$ is reachable by $(\overline x, y)$  
   from  $s$  
   by a    move $b$  of Player Two with $g^b<\overline \ep$:}

We have $p^s = 1/Q_2$, $p^s\esc (u,s)\leq g^b \leq 2.3 a(s)$ and  
$${\esc(u,s) \over \mu(u,s)} \leq {2.3 \ a(s) Q_2 \over \mu (u,s)} 
\leq  {2.3 \ a(s) Q_2 \over L \tilde   a (s)} 
\leq {2.3 \ \over Q_1}.$$

In all arguments that follow concerning 
 members of a set $T$ as created above, for 
 convenience we will write $z^{\overline \ep}$ or $B^{\overline \ep}$ 
instead of 
 $z^{\gamma}$ or $B^{\gamma}$ for $\gamma >\overline \ep$.
 By Lemma 5.5 there will be no difference 
  in these expressions.  \vskip.2cm

{\bf Step 2; For any choice of $s\in R$ from Step 1
there is an $\tilde \ep$  simplication 
 $\overline y$ of Player Two's strategy $y$ such that 
together with $\overline x$ the state $s$ and all states $t\in T$ with 
 $z^{\overline \ep}_{\overline x, y}(t)
 \geq \hat \ep \ \overline \ep \  z^{\overline \ep}_{\overline x,y}(s) 
 /(4|{\cal N}|)$ 
 are reached by  $(\overline x, \overline y)$ from all of $R$,
  and furthermore from inside of $T$ no state 
 outside of $T$ is reached:}  \vskip.2cm

We define $\overline y^t$ for all $t\in T$
 by removing from $y^t$ all moves made by 
 Player Two with a frequency of  $L/ (L-1)Q_1$ or less, followed 
 by normalization.  

  Any $t\in T$  
     that satisfies   
 $z^{\overline \ep}_{\overline x,y} (t) \geq  
   { \hat \ep \ \overline\ep  \ \overline z^{\overline \ep}_{\overline x,y} (s)
 \over  4 |{\cal N}|}$
 by Step 1    also satisfies 
    $P^s_{\overline x,y} (s,t) \geq 
  { \overline \ep \ \hat \ep \ \beta \over 4.5|{\cal N}|}
    P^s_{ x,y} (s,t)$ and
    $z^{\overline \ep}_{x,y}(t)\geq 
{\overline \ep\ \hat \ep \ \beta  \over 4.5|{\cal N}|}$.
    Notice that this last condition is  
      satisfied by the state $t=s$.
 For any $t\in T$ with $P^s_{\overline x,y} (s,t) \geq 
  { \overline \ep \ \hat \ep \ \beta \over 4.5|{\cal N}|}
    P^s_{ x,y} (s,t)$ and  $z^{\overline \ep}_{x,y}(t)
\geq {\overline\ep \ \hat \ep \ \beta \over 4.5|{\cal N}|}$ to 
 show that $t$ is reached from all of $T$ with $(\overline x , \overline y)$
 by Lemma 3.4 it suffices 
   to show that for any $w\in T$ and any $t\in T$ satisfying 
$z^{\overline \ep}_{x,y}(t) \geq {\overline \ep \ \hat \ep\  \beta \over 
 4.5|{\cal N}|}$, including $s=t$, we have that the change from 
 $y^w$ to $\overline y^w$ doesnot reduce 
 $P^w_{x,y}(w,t)$ by more than a factor of $\overline \ep \hat \ep \beta / 
(12 |{\cal N}|^2)$.

      If $b\in A^w_2$  
  is a Player Two move with $g^b\geq \overline \ep$,
 removing $b$ to form $\overline y^w$ from 
 $y^w$ cannot 
       reduce     $P_{ x,y}^w(w,t)$ by anything more than a factor of 
 $\ep^*/\overline \ep$. 
       Assuming that $g^b<\overline \ep$
 and removing $b$ to make 
 $\overline y^w$ removes 
        at least 
        ${\overline\ep \ \hat\ep \  \beta \over 12 m|{\cal N}|^2}$
         of the motion 
         $P^w_{x,y} (w,t)$ 
         we would have from Lemma 5.5  
         that $2.3 \   a(w) L 
         /Q_1 (L-1) \geq g^b_{x,y}L/(L-1)Q_1 \geq g^b_{x,y} y_b\geq     
 {\overline\ep \ \hat\ep \ \beta \over 12 m|{\cal N}|^2}
 {\overline \ep \ \hat \ep \  \beta  \over 4.5 |{\cal N}|}
 a(w)$.  
         This is a contradiction to the definition of $Q_1$.

Second, we show that, starting at $s$, motion according  
 to $(\overline x, \overline y)$ never leaves the set $T$. 
    Let us assume  that 
 $u$ is a state not in  $ T$  
reached by a  move $b$ of 
 Player Two from any  $t\in T$     played      
   against $\overline x$ and given positive frequency by $y$. 
              We need to show that 
   $b$ is not used in $\overline y$.
 If $t \neq s$ then                          by formula (3)
$ \mu(t,u) \leq {a(t) \over p^t y_b} = {\mu(s,t) Q_2 \over y_b}. $
In particular,  by the definition of $T$ and since  $\mu$ is a metric, 
$$y_b \leq { \mu(s,t) Q_2 \over \mu(t,u)} \leq  
{ \mu(s,t) Q_2 \over \mu(s,u)} { \mu(s,u)  \over  \mu(t,u)}\leq 
{Q_2 \over L-1} = {L \over (L-1) Q_1}.$$
And if $b$ is a move  at the state $s$ then also by the definition 
 of $T$                      and (3)
$$y_b \leq { a(s) Q_2\over \mu (s,u)}
\leq {Q_2 a(s) \over L \tilde   a(s)} 
\leq {1\over Q_1}.  $$ 
       Therefore we conclude that $(\overline x, \overline y)$ defines 
        one ergodic set $D\subseteq 
 T$ that includes $s$ and all states $u\in R$ 
         satisfying 
        $z^{\overline \ep}_{\overline x ,y}(u) \geq 
       {\overline\ep  \ \hat\ep \  \over 4|{\cal N}|}$.  
\vskip.2cm 

{\bf Step 3; Show that there is a proper choice of $s$ from Steps 1 and 2 
 with a subset $V_R$ of Player Two moves satisfying 
  the conditions of Theorem 1, namely that these moves belong to a subset $F$  
 containing $s$ and inside 
of the ergodic set $D$  such that $\tilde w (t)\tilde a(t)$ is 
 a constant for all $t\in F$ and there is a distribution on $V_R$ such 
 that  used
 against $\overline x$ gives an expected payoff to Player One within 
 $ \hat \ep$ of $r^1(s)$:} 
\vskip.2cm 

Define $U:= \{ t\in R\ | \ \xi (t)\leq r^2(t)-2.4 
\overline \ep\}$ and define $\tilde U:= \{ t\in R\ | \ \xi(t) \leq 
 r^2(t)-2.5 \overline\ep\}\cap \{ t\in R \ | \ z^{\overline \ep}_{x,y}(t) 
 \geq \beta \}$. 
 We create a partition $\{ U_1, \dots ,U_k\}$ of the members of 
 $ U$ in increasing values of $\tilde w \tilde a$,
 meaning   
 that $s$ and $t$ belong to the same  
  member of  $U_i$ if and only if $\tilde w (s) \tilde a(s)
=\tilde w(t) \tilde a(t)$. 
    For any state $s$ in $\tilde U$  
     we consider the sets $T(s)$ and $D(s)$ and the  strategies 
    $\overline x (s), \overline y(s)\in X\times Y$ 
 as created above in Step 1 and Step 2.

For the sake of contradiction we suppose that there is no  $s\in  
 \tilde U\cap U_i $
 and  
$b\in B_{x,y}^{\overline \ep}(t)$ with $t\in D(s) \cap U_i$
 such that $|v^1_{\overline x(s),y}(b)-r^1(s)|\leq 
 \hat \ep $ and 
 there is no pair of Player Two moves $b,b'\in B^{\overline \ep}_{x,y}(R)$
 with both  $b$  and 
 $b'$ belonging 
 to the  set $D(s)\cap U_i$ with $v^1_{\overline x(s),y}(b)$
 and $v^1_{\overline x(s),y}(b')$ on 
 different sides of $r^1(s)$.

For every $s\in \tilde U$ and 
$t\in U\cap D(s)$ with some 
 move in $B_{x,y}^{\overline \ep}(t)$ used in $y^t$  let  $v_s^1(t)=
\sum_{b\in B^{\overline \ep}_{x,y}(t)} v_{\overline x(s),y}^1(b) 
\nu^b_{\overline x(s),y} /$ 
$\sum_{b\in B^{\overline \ep}_{x,y}(t)}  \nu^b_{\overline x(s),y} $, 
 the average Player One payoff resulting from these moves at $t$. 
 For every $1\leq i\leq k$ let $p(i):=\sum_{j<i}|U_j|$.

We claim that our above assumption implies that 
  $ z^{\overline \ep}_{x,y}(s)
\leq 3^{p(i)}\beta /( \hat \ep)^{p(i)}$
 for every $s\in U_i\cap \tilde U$.  

We prove the above claim by induction on $i$. 
Let $s$ be any member of $U_i\cap \tilde U$, and we  
 assume that $|v_s^1(s)-r^1_{\overline x(s),y}(s)|=
 |v_s^1(s)-r^1_{ x,y}(s)|\geq \hat \ep$.
From Part 1  and Part 2 we know that 
 the importance with respect to $(\overline x(s),y)$ 
 from exits  outside of $B^{\overline \ep}_{x,y} (D(s))$ 
 does not exceed $\overline \ep \ \hat \ep/3$ times 
 the importance of the $B^{\overline \ep}_{x,y}(s)$ moves. 
 Since all the $v_s^1(t)$ 
 with $t\in D(s)\cap U_i$ 
 are on the same side of $r^1(s)$ as $v_s^1(s)$, we are left 
 only with the $B^{\overline \ep}_{x,y}$ moves 
  from $\cup_{k<i}U_k\cap D(s)$ to counter-ballance the 
 $v_s^1(s)$ to make $r^1_{\overline x(s),y}=r^1_{x,y}$. 
We can assume now 
 that $i>1$, since otherwise we would have to conclude that 
 $|v_s^1(s)-r^1_{ x,y}(s)|\geq \hat \ep$ is impossible. 
By the induction hypothesis 
 the sum of all the $z^{\overline \ep}_{x,y}(u)$ over the set 
 $ \cup _{k<i} U_k$ does not exceed    
 $\sum _{k<i}|U_{k}| 3^{p(k)}\beta / \hat \ep^{p(k)}\leq 
 {2\over 3} 3^{p(i)} \beta / \hat\ep ^{p(i-1)}$. By the fact 
that our simplications $\overline x(s)$ hardly influence 
 the expected payoffs from moves with an absorption rate 
 of at least $2\overline \ep$ and by the statement 
 of Step 1,  in order  
 to maintain $|v_s^1(s)-r^1_{ x,y}(s)|\geq \hat \ep$ we must 
 assume that $\hat \ep z^{\overline \ep}_{x,y}(s)\leq 
 {3\over 2} {2\over 3} 3^{p(i)} \beta / \hat\ep ^{p(i-1)}$, and this 
 concludes the proof of our claim.

With the definition 
 of $\beta$ we conclude 
 that $z_{x,y}^{2.5 \overline \ep}(s)<\tilde \ep /|R|$ for 
 every $s\in R$, and this means that $R$ could not have been chosen 
 for polarization, a contradiction.  
  \vskip.2cm

 With the appropriate $s\in R$ chosen, we have $D_R:= D(s)$, 
 $x_C $ defined from $\overline x(s)$ and  
 $y_C$ defined from $ \overline y (s)$ so that changes are made 
 only   inside of $D_R$, and the exits 
 $V_R$ and their distribution as determined by $y_D$ 
 come from the above argument. 
\vskip.2cm 

{\bf Step 4; show that  the moves
  $B_R$  satisfy the requirements of Theorem 1:}
 \vskip.2cm 
 
The easiest way to prove that 
 $|r_{x,y_R}(s_R)-r_{x,y}(s_R)|\leq \hat \ep$ is to return 
 part of the way back to the space ${\cal S}_{\sharp}$!  We let 
 $\tilde {\cal S}_{\sharp}$ be the space generated by the almost 
 trivial partition 
 $\tilde {\cal P}:= \{ R\} \cup \{ \{ s\} \ | \ s\not\in R\}$. With 
 $\tilde r^1$ the  harmonic function on $\tilde {\cal S}_{\sharp}$ 
 induced by $r^1$ on the absorbing states,
 by Lemma 3.5  $\tilde r^1(s_R)$ and  $r^1(s_R)$ differ 
 by at most 
 $4\ep ^*$. Let $\tilde \nu_{\sharp}$ be the corresponding 
 measure of the importance of the exits.

 Define a move $a\in A^s_1$ 
of Player One in the set $R$ to be a {\em principle} move if $a$ is not 
 a limit move and if
 there is a $b\in   B^s_2$ such that $(a,b)$ is an exit with
 $\nu_{x,y}(a,b)\geq  \hat \nu - \ep^*$.

We claim 
 that the combination $(a,b)$ of a move of $B_R$ with a principle move 
 of Player One must yield  $\nu(a,b)\leq \ep^*$.  Once this 
 is established from the definition of $B_R$ 
 we need only to break down   the sum of 
  the $v^{\tilde r^1}(a,b) \tilde \nu_{\sharp}(a,b)$ over 
 all exits $(a,b)$ with $\nu(a,b)\geq \hat\nu -\ep^*$ 
and apply Lemma 3.7 to conclude that   
  $r^1_{x,y_R}(s_R)$ is within $20|{\cal N}| m \ep^*/\hat\nu$ 
 of $r^1_{x,y}(s_R)$,  that is  
 much closer than we need it.  Suppose for the sake 
 of  contradiction that for some principle $a$ and some $b\in B_R$ that 
 $\nu (a,b)\geq \hat \nu-\ep^*$.  Assuming that the moves 
 take place at $t$, we have 
from the definition of $B_R$ that $a(t)\geq y_b (\hat q-\ep^*) (\hat \mu - 
\ep^*) \rho$.  Furthermore by definition 
 we have $\nu(a,b) \leq x_a y_b / a(t)$ and by assumption 
 $x_a\leq \ep^*$. These four inequalities 
 are contradictory.

We  show that $b\in B_R\cap 
 A^t_2$ with $t\in P$ implies  
 $v^2(b)< r^2(t)$. If $\xi(t)\leq r^2(t)-2\overline \ep$ then 
it follows from Lemma 5.5.  If $\xi(t) > r^2(t)-2\overline \ep$ 
 then by Lemma 5.4 all moves have the same auxillary value 
 $\xi (t)= \overline \xi (t)$; 
it follows  from the smallness 
 of $\delta\leq  (\ep^*)^2$ and formulas (10) and  
 (11) that if  $v^2(b)\geq r^2(t)$ then the repeated use of $b$ 
 would result in a higher evalution for $\xi(t)$ because an undiscounted
 value of at least $r^2(t)$ would be obtained but at much higher
 auxillery absorbing rate.
  \vskip.2cm

{\bf Step 5; show $z^{2.5\overline \ep} (s)< \tilde \ep$ for any 
 state $s$ that is not in $P$ or is a satelite of some $R\in {\cal R}$:} 
\vskip.2cm  

If $s$ is not a satelite and not in $P$ then 
 due to the very small size of $\delta$ we have 
 from (11) that  $\xi(s)$ is within $\tilde  \ep$ of 
 $r^2(s)$, implying that no move $b$ used at $s$ could 
 satisfy $w^2(b)< r^2 (s)-2\overline \ep$.   
For a satelite $s$ of $R$ we suppose that  
  $b\in A^s_2$ is a Player Two move at $s$
  with $g^b\geq 2.5\overline \ep$.  Such moves have 
  at least a $2.4\overline \ep$ probability of never 
 returning to the set $R$.  The probability of using such a 
 move before reaching $R$ must be no more than $\ep^*/(2.4 \overline \ep)$, 
and thus the total probability that it is used cannot exceed 
 $\ep^*/ (2.4 \overline \ep (\hat \mu -\ep^*))$. 
    \hfill q.e.d.\vskip.2cm

\section {Signaling}

In this section we show that there are approximate 
 equilibria without an assumption that  Player One can send 
  signals independent of the transitions. 
 The problem concerns the consequences to the players of any moves that 
 would be used by Player One as a transition dependent signal.  For 
 example, a move of Player One  that 
 brings the play outside of the set $D_R$ may fail to be useful to signal 
 her desire for Player Two to use a move in $V_R$, because outside
 of $D_R$ the jump function for Player Two may exceed greatly 
 his expected payoff 
 from the moves in $V_R$.

The natural solution is for Player One to have a  move   
 inside of $D_R$ that is not used in $x_C$ whose use means    
 that the  moves $V_R$ 
 of Player Two will not be used, and after a certain 
 quantity of visits to some state in $D_R$   
 it will be understood mutually  
 that Player Two must use a move in $V_R$.   
 A problem arises, however, if every such move results in a positive  
 probability of leaving the set $R$.

With regard to the next two theorems, we assume  the statement and proof of 
 Theorem 2, which means also that we assume that 
all the conditions of Theorem 1 are satisfied. 
 We will add new conditions to those of Theorem 1 and 
   make some  minor changes to  
  the proof of Theorem 1. The definition of ${\cal S}_{\sharp}$ remains, along 
with its  Markov chain transitions, including the $p^*_R$ and $p_R$. 
 The changes begin with the definition of the parts 
 $q^d_R$ and $q_R$  and therefore everything that follows in 
 the proof of Theorem 1 will be altered as well, including 
 the introduction  of   new situations. 
\vskip.2cm 

{\bf Theorem 1':} Assume the following 
 property for every $R\in {\cal R}$: if 
  every  move $a\in A_1^t$ in $ D_R$ removed to make $ x_C$ from $x$ 
 formed an exit against some Player Two move used in $ y_C$, 
 then there exists a set $A_R$ of Player One principle moves in 
 $D_R$ such that  \newline 
1) the sum of $\nu_{x,y} (a,b)$ for all $R$  exits $(a,b)$ performed 
 outside of $D_R$   
   does not exceed 
 $\tilde \ep \ \overline \ep \  \hat\ep\ \beta /3$, \newline 
2) for every principle move $a\in A_R$ of Player One                        
 used at $t\in D_R$            with 
$\nu^a _{x,y} \geq \beta \ \hat \ep \ \tilde \ep\ \overline \ep / 
 (3|{\cal N}| m)$ we have 
$\sum_{\ b \mbox { used in }y^t_C} 
\nu_{x,y}(a,b)\geq (1-\overline \ep \ \hat \ep \  \tilde \ep\ \beta)
 \nu^a_{x,y}$ and therefore also 
$|v^1_{ x_C ,  y_C}(a)-v^1_{ x,y}(t)|\leq 
   \hat\ep$.\newline 
{\bf Conclusion:} Without any assumption 
 on Player One's ability to signal indendependently of 
 the transitions,  the game has approximate 
 equilibria.     
\vskip.2cm

{\bf  Proof:}
Define a member of ${\cal R}$ to be {\em problematic} if 
 the assumption of Theorem 1' holds. 
We proceed exactly as the proof of Theorem 1, except that for 
 all problematic $R$  
 we incorporate into the ${\cal S}_{\sharp}$ transition $q^c_R$
 all the  $R$ exits not inside of $D_R$   
  or not created from a combination    
 of an  $a\in A_R$ with a move used in $y_C$. 
Recalling that $q^d_R$ is the difference between 
 $q^c_R$ and $p^*_R$ by    Lemma 3.7  we still have that 
  $\nu_{\sharp} (q^d_R) (v^2_{\sharp}(q^d_R)- r^2_{\sharp}(s_R))$ is 
 below $\hat\ep $. Due to  
 Condition 2 and Lemma 3.7 we have the other requirement for 
 applying Lemma 3.9. 
 We assume that  ${\cal T}$ is the subset of ${\cal R}$ that has been 
 polarized.

  Define a situation $s^w$ at a state $s$ to be {\em timed} if there is a 
 natural number $m$   
  such that $s^w$ is determined by the present state $s$ and 
 the previous situations and moves in  
 the last $m$ stages.  A normal situation is timed, but the converse doesn't 
 hold.  

We keep  the same situations $s^e$, $s^f$ and $s^g$ from   
the proof of Theorem 1.  The stationary strategies 
 for all the $s^g$ and all 
the $s^e$ other than a representative  $s^e_R$ are defined in the  
 same way, and in a non-problematic $R$ the stationary strategies for 
 $s^f$ are also the same.   

 For every polarized  $R\in {\cal T}$ and $t\in D_R$ we create 
 a timed situation $t^h$.  
When  a situation $s^e_R$ is reached the strategies $(x_C, y_C)$ are
 performed, 
 but   instead of moving to 
 a $t^f$ or $t^g$ there is motion to the timed situation $t^h$.

For non-problematic polarized 
 $R\in {\cal T}$ we choose any $t\in D_R$ such that 
 there is a Player One move $a$ at $t$  not used in $x_C$ 
 and when  paired with $y_C$ results in zero probability of leaving 
 the set $R$. Create a frequency  $\tilde f_a>0$ and a number 
 $n_t$ such that $f_a\sum_{i=0}^{n_t-1} (1-f_a)^i= 1-\lambda_R$, where 
 $\lambda_R$ is that quantity determined by the polarization, and 
 such that for any distinct $u,v\in D_R$ the probability of using 
 the move $a$ before  moving 
 from $u$ to $v$  is at least $1-\ep^*$.    
For all the situations $s^h$ for $s\not=t$  the players act according 
 to $(x_C, y_C)$  
 and at $t$ Player Two  according to $y_C$ and 
  Player One according to $(1-f_a) x_C+ f_a 1_a$.  
 If on the $n_t$th visit to the situation  $t^h$ the move $a$ was not made, 
then the situation following $t^h$ is some $u^g$.  Otherwise 
 if the move $a$ was used on any visit to the situation 
 $t^h$ then the next situation is either some    
     $u^f$ if an exit wasn't used or some $u^e$ if an exit was used.

For problematic $R\in {\cal R}$, 
 let $\pi_R\in \Delta (A_R)$ be 
 the probability distribution on the $A_R$ 
  that is  generated conditionally by 
 $(x,y_C)$.
  Choose a natural number $n_R$ and 
a stationary strategy $x_C^*$ for Player One 
 so that with a start at $s_R$ the distribution on the moves $A_R$ 
 is $\pi_R$  and for every pair $u,v\in D_R$  the probability of
 using a move in $A_R$ before moving from $u$ to $v$ is 
 no more than $\ep ^*$ and the probability of using 
 some member of $A_R$ at or before 
 the $n_R$th visit to the state $s_R$ is $1-\lambda_R$.       
 For the situations $t^h$ with $t\in D_R$ the players act according 
 to $(x^*_C, y_C)$.  
 If on the $n_R$th visit to the situation  $s_R^h$ the move $a$ was not made, 
then the situation following $t^h$ is some $u^g$.  Otherwise 
 if a move in $A_R$ was used on any visit to a situation 
 $t^h$ then  the next situation is either   
     $u^f$ (if an exit was not used) or $u^e$ (if an exit was used). 
   At a situation $u^f$ the strategies $(x^*_C,y_C)$ are also used.  

As with the  proof of Theorem 1 we must  show that
 the expected 
 payoffs to Player $i$ from  
 every  situation $s^e$ is within $3.1\overline \ep$        
of $r^i_{x,y}(s)$. 
Given the proof of Theorem 1 the only additional 
 argument needed  concerns the 
   use of  exits in a problematic $R$ 
  before  the timed situations have been 
 reached. This did not present a problem in the proof of Theorem 1 
 because they were the same exits used in the situations $t^f$ and 
 performed with the same distributions. If we can show that the 
 total probability of their occurance cannot exceed $\overline \ep /10$, 
 then we get our result by  ignoring their influence.    
Indeed in the Markov chain defined on ${\cal S}_{\sharp}$ 
 the absorption rate of $s_R$ for a problematic $R$ 
 is at least $ \rho \hat \mu/(2Q_1) $.
  By Lemma 3.9 this absorption rate 
 does not fall below ${\rho \hat \mu \over 2Q_1}
{\overline \ep ^{3|{\cal N}|} \over (3|{\cal N}|)^ 
 {|{\cal N}|}}$ after polarization. Since this quantity is 
 still very large compared to $\ep^*$, the maximal probability 
 of using such a exit before a timed situation is reached, 
 we can indeed ignore these exits.
  (We leave the formal argument using  Section 3 to 
 the reader.)

 The situations  defined  above are not normal and thus 
 do  not generate a stochastic game, preventing 
 a direct application of Corollaries 4.3 and  4.4. 
 Therefore we must  perceive the situations 
 $\{ s^h\ | \ s\in R\}$ for $R\in {\cal T}$  as  subgames. 
Concerning the behavior of Player One, we  view the entire 
 process up until the $n_R$th visit to the state $s_R$ or the 
 $n_t$th visit to $t$   
  as a single 
 decision -- whether or not to use a  move in $A_R$ and  
  if so then 
 which one. This places Player One's decisions back into 
 the context of Corollary 4.4.  
  
 Concerning the behavior of Player Two, the matter is  
 more complex.  Player Two could have an  influence on the payoffs 
 by altering  the strategy $ y_C$. 
 Strictly speaking the context would be no longer that of a 
 harmonic function  on a time 
 homogeneous Markov chain  -- the expected payoff to 
 Player Two at a   state corresponding to a situation 
 $t^h$  would be  changing over time. 
However Player Two's ability to gain or lose in expected payoff 
 is  conditioned on the  use of a move of Player One in 
 $A_R$ --  this  {\em is} modeled by a time homogeneous 
 Markov chain and therefore   Proposition 4.2 is sufficient for 
 the conclusion. 
       \hfill $\Box$\vskip.2cm

{\bf Theorem 2'} 
 The conditions of  Theorem 1' are  satisfied always.
  \vskip.2cm

{\bf Proof:}
Let $(x,y)\in X\times Y$   
 and $(x_C, y_C)$ be a solution 
 given by Parts 1, 2, and 3 of Theorem 2 for a 
  polarized $R\in {\cal R}$ and we assume that 
 conditions of Theorem 1' hold for $R$ (meaning that 
 $R$ is problematic).  

{\bf  1)} 
 Consider  the strategies played at any $t\in D_R$.    
Suppose for the sake of contradiction that there is a state 
 $u\in R\backslash D_R$ where an importance of at least 
$\beta \ \tilde \ep \ \overline \ep \ \hat \ep /3|R| $  
   occurs from exits at $u$. 
Consider the moves that were 
 removed from $y^t$ to make $\overline y^t$.   
   By  Lemma 5.5  
        at any $t\in D_R$ no more than 
       $ { 7|R|\over \tilde \ep \overline \ep \ \hat \ep\ \beta}
\ {mL\over (L-1)Q_1}$ of 
       the transition $P^t(t,u)$ was removed to make $ y_C^t$ from  
        $y^t$.  On the other hand, given that every   
 move of Player One removed from $x^t$ to make $ x_C^t$  
 would have created an exit against some move in $ y_C^t$, 
   we must also conclude from the rare use of an exit  
 that no more 
 than       $2 \ep^* Q_1$ 
        of the transition  in $P_{x,y}^t(t,u)$ came from 
 such a Player One move.  From Lemma 3.3 we have  
        that $u$ is in $D_R$, a contradiction. 
 
{\bf 2)}
Assume that $\nu_{x,y}^a \geq \beta\ \hat \ep 
 \ \overline \ep \ \tilde \ep / (3|{\cal N}|m)$ 
 for some principle move $a$ of Player One at $t\in D_R$.
Suppose for the sake of contradiction   
 that the probability of reaching any absorbing 
 state  
 from this principle move is altered by a 
factor of more than $\beta \overline  \ep  \tilde \ep \hat \ep /2 $ 
  by the change from $y$ to $y_C$. This means that  
$\nu _{x,y} (a,b)$   
 is at least  $\overline \ep ^2\ \tilde \ep ^2 \ \hat
 \ep^2 \ \beta ^2 \over 6|{\cal N}|^2 m^2$ for at least one 
 move $b$ that was removed to make $y^t_C$ from $y^t$. We must 
 conclude from Lemma 5.5 that ${\beta^2 \ \hat \ep ^2\ \overline \ep^2 \
 \tilde \ep^2 
   a(t)\over 6|{\cal N}|^2m^2}\leq {2.3 a(t)L \over Q_1 (L-1)}$,
 a contradiction to the definition of $Q_1$.
   
 The final claim follows now from the  
 argument in part 4 of the proof of Theorem 2, showing that 
 $v^1_{x,y}(a)$ is very close to the value of $r^1$ for all primary 
 moves. 
  \hfill $\Box$
\vskip.2cm 

In the proof of Theorem 1' we could eliminate the argument that exits 
 performed before reaching a timed situation in a problematic set 
 are irrelevant if   
 we had  a more powerful   Markov chain result (that combines 
 the condition of Lemma 3.3 with the conclusion of Lemma 3.2)    
   or we use  Vieille's 
approach  to  ``communication sets" (Vieille 
 2000c),  showing how one can move through a set $R$ with no danger 
 of leaving it.    

\section{Countably many states}

On the  technical side, the problem with applying either our or 
 Vieille's proof to countably many states lies in the 
 finite state space assumption that given 
 any stationary strategies for the players and  any positive 
 $\delta$ there will exist  
 a $\delta$ perturbation of this strategy that is absorbing.  

A strategy for finding a counter-example could be following.  Construct 
 an infinite sequence of   
 games $\Gamma_0, \Gamma_1, \dots$ that are positive recursive 
 for both players
 corresponding to  increasing finite sets $S_0\subseteq S_1\subseteq \dots $ 
of  non-absorbing states 
 such that  for every $i\geq 0$ and $j\geq i$ the moves and their induced 
 motions inside of $S_i$ are the same for all games $\Gamma_j$. 
Construct a countable state space by having the game start at $s_0$, 
 define the state space on the $i$th stage to be the space $S_i$, and 
 declare that absorption occurs on stage $i$ if an absorbing state 
 of the game $\Gamma_i$ is reached.  Furthermore,
 give both players the  
 ability to force the game to absorption in the new countable 
  state space game.  
Desirable may be  games $\Gamma_i$ such that 
 with large $i$   the approximate  equilibrium behavior 
 of $\Gamma_i$  keeps the non-absorbing 
play  most of the time 
 close to  the set  
   $ S_0$ and 
 the minimal number of stages necessary to reach an absorbing 
 state in the game $\Gamma_i$ starting 
 from any  $s_0\in S_0$ goes to infinity as $i$ goes to infinity.  Otherwise 
 if we allow that absorbing states are reachable quickly 
 from all non-absorbing states, 
to avoid convergence toward  large sub-games of  
essentially equilivalent states it  may be desirable  
  if reaching an absorbing 
 state of $\Gamma _i$ on the $i$th stage of play does not mean 
 certain absorption but rather a positive probability of  absorption   
  mixed with a positive  
  probability  of starting  the game over   
 at  $s_0\in \Gamma_0$.   

There are many ways for a game to have a countable state space but 
 be played essentially on finitely many situations, for example games 
 that break down into sequences of sub-games played essentially 
 on finite state spaces.  
 Also to be avoided are structures that are formally countable in size but do 
 not exploit the full potential of what it means to have infinitely many 
positional possibilities.  We believe that the best 
 candidates for a counter-example    will 
 incorporate the concept of a random walk on arbitrarily 
 many positions, as presented  
 in  our introduction. However, 
  to avoid  operator approaches similar to that of the Maitra and Sudderth 
  proof   we believe 
 that there  must be a conflict by {\em both} players 
 between exploiting their positions and controlling   the behavior of the 
 other player.  For this and other reasons, we believe 
that  the non-absorbing states   
 must have a structure more complex than ${\bf Z}$, for 
 example involving joint random walks  on the   
   two dimensional lattice 
 ${\bf Z}^2$.

\vskip1cm 
{\bf Additional Acknowledgment:} 
 The  author  thanks 
   Cafe Europe on Azzahra Street in Jerusalem for the many hours       
    he worked on the proof at the cafe.

\section {References}
\begin{description}
 \item [Bewley, T. and Kohlberg, E. (1976).]
   ``The Asymptotic Theory of Stochastic Games",
   {\it Math. Oper. Res.},
   {\bf 1}, pp. 197-208.
\item [Blackwell, D.  and Ferguson, T (1968).] ``The Big Match", 
{\em Annals of Mathematical Statistics}, {\bf 39}, pp. 159-163.
\item [Fink, A.M. (1964).]
   ``Equilibrium in a Stochastic $n$-Person Game",
   {\it J. Sci. Hiroshima Univ.}, {\bf 28}, pp. 89-93.
\item [Freidlin, M. and Wentzell, A (1984).] {\em Random Perturbations 
 of Dynamical Systems,} Springer, Berlin.
\item [Kohlberg, E. (1974).]
   ``Repeated Games with Absorbing States",
   {\it Ann. Stat.}, {\bf 2}, pp. 724-738.
\item [Maitra, A. and Sudderth, W. (1991).] ``An Operator Solution of 
 Stochastic Games", {\em Israel Journal of Mathematics}, {\bf 78}, pp. 33-49. 
\item [Martin, D. (1975).] ``Borel Determinacy", {\em Annals of 
 Mathematics}, {\bf 102}, pp. 363-371.
\item [Martin, D. (1998).] ``The Determinacy of Blackwell Games", 
 {\em The Journal of Symbolic Logic,} {\bf  63}, pp. 1565-1581.
\item [Mertens, J.F. (1987).] ``Repeated Games", in {\em Proceedings of the 
 International Congress of Mathematicians}, (Berkeley, 86),
 American Mathematical Society, pp. 1528-1577. 
\item [Mertens, J.F. and Neyman, A. (1981).]
   ``Stochastic Games",
   {\it Int. J. Game Th.}, {\bf 10}, pp. 53-66.
\item [Shapley, L.S. (1953).]
   ``Stochastic Games",
   {\it Proc. Nat. Acad. Sci. U.S.A.}, {\bf 39}, pp. 1095-1100.
\item [Solan, E. (1999).] ``Three-Player Absorbing Games", 
 {\em Mathematics of Operations Research}, {\bf 24}, pp.  669-698. 
\item [Solan, E. (2000).]
   ``Stochastic Games with 2 Non-Absorbing States",
   {\em Israel Journal of Mathematics}, {\bf 119}, pp.  29-54.
\item[ Thuijsman, F. and Vrieze, O.J. (1989).]
   ``On Equilibria in Repeated Games With Absorbing States",
   {\it Int. J. Game Th.}, {\bf 18}, pp. 293-310.
\item [Vieille, N. (2000a).]
   ``Two-Player Stochastic Games I: A Reduction",
   {\it Israel Journal of Mathematics}, {\bf 119}, pp. 55-91.
\item [Vieille, N. (2000b).]
   ``Two-Player Stochastic Games II: The Case of Recursive Games",
{\it Israel Journal of Mathematics}, {\bf  119}, pp. 92-126.
\item [Vieille, N. (2000c).]
   ``Small Perturbations and Stochastic Games",
   {\it  Israel Journal of Mathematics}, {\bf 119}, pp. 127-142.
\item [Williams, D. (1991).] {\em Probability with Martingales}, 
 Cambridge University Press.
\end{description}

\end{document}